\documentclass[12pt]{article}
\usepackage{amsmath}
\usepackage{amsfonts}
\usepackage{amssymb}
\usepackage[latin1]{inputenc}

\title{Irrationalit\'e de valeurs de z\^eta
[d'apr\`es Ap\'ery, Rivoal, ...]}
\author{St\'ephane FISCHLER}
\date{S\'eminaire Bourbaki - Novembre 2002 \\ 
Expos\'e num\'ero 910 ; \`a para\^{\i}tre dans Ast\'erisque}

\def\@pnumwidth{1.6em}
\def\@tocrmarg {2.6em}
\newcounter{secuppdepth} 
\newtheorem{defi}{\definame}[section]
\newtheorem{prop}[defi]{\propname}
\newtheorem{theo}[defi]{\theoname}
\newtheorem{conj}[defi]{\conjname}

\newtheorem{lemm}[defi]{\lemmname}
\newtheorem{rema}[defi]{\remaname}

  \def\@byname{par}

  \def\definame{D\'efinition}
  \def\propname{Proposition}
  \def\theoname{Th\'eor\`eme}
  \def\conjname{Conjecture}
  \def\coroname{Corollaire}
  \def\faitname{Fait}
  \def\exemname{Exemple}
  \def\lemmname{Lemme}
  \def\remaname{Remarque}
  \expandafter\ifx\csname fhyph\endcsname\relax
    \else\fi 

\newcommand{\Pun}{{\mathbb P}^1}
\newcommand{\moins}{\setminus}

\newcommand{\eps}{\varepsilon}
\newcommand{\croix}{\times}

\newcommand{\aeno}{\ell}

\newcommand{\fpqusz}[9]{{_{#1}F_{#2} \left(
\begin{array}{ccccc|}
 #3, & #4, & #5, & \ldots, & #6 \\
    & #7 ,& #8, & \ldots, & #9 
\end{array}
\, \, \,  z^{-1} \right)}}

\newcommand{\fquatretrois}[1]{{_4 F_3 \left(
\begin{array}{cccc|}
 -n , & -n,  & n+1, & n+1 \\
     &  1 , &  1 , &  1 
\end{array}
\, \,  #1 \right)}}

\newcommand{\fquatretroispre}[1]{{_4 F_3 \left(
\begin{array}{cccc|}
 -n , & -\yg, & n+1, & \yg+1 \\
     &  1,  &  1, &  1 
\end{array}
\, \,  #1 \right)}}

\newcommand{\ope}{L}  
\newcommand{\rivo}{{\bf L}} 
\newcommand{\Lmod}{{\mathfrak{D}}}  
\newcommand{\de}{\delta}
\newcommand{\lam}{\lambda}
\newcommand{\gam}{\gamma}
\newcommand{\om}{\omega}

\newcommand{\AAA}{{\bf B}}
\newcommand{\BBB}{{\bf A}}

\newcommand{\RRR}{{\bf R}}  
\newcommand{\SSS}{{\bf S}}  
\newcommand{\Rvetun}{{\bf \overline{R}}_n} 
\newcommand{\Svetun}{{\bf \overline{S}}_n} 
\newcommand{\Svetunphi}{{\bf \overline{S}}_{\varphi(n)}}  
\newcommand{\cijvetun}{\overline{c}_{i,j}}   
\newcommand{\Pvetun}{\overline{P}}  
\newcommand{\Rzu}{\widetilde{{\bf R}}_n}
\newcommand{\Szu}{\widetilde{{\bf S}}_n}
\newcommand{\cijzu}{\tilde{c}_{i,j}}
\newcommand{\Pzu}{\widetilde{P}}
\newcommand{\BB}{H}
\newcommand{\QQ}{Q}
\newcommand{\Fp}{{\mathbb F}_p}
\newcommand{\Gm}{{\mathbb G}_m}
\newcommand{\Asoro}{T}
\newcommand{\Bsoro}{U}
\newcommand{\Csoro}{V}
\newcommand{\Dsoro}{W}
\newcommand{\Aorth}{\widetilde{A}}
\newcommand{\Borth}{\widetilde{B}}
\newcommand{\Corth}{\widetilde{C}}

\newcommand{\pre}{\varphi}
\newcommand{\yg}{x}
\newcommand{\Ppre}{P}
\newcommand{\Qpre}{Q}

\newcommand{\cij}{c_{i,j}}
\newcommand{\gdo}{\textup{O}}
\newcommand{\dd}{{\rm d}}
\newcommand{\Li}{{\rm Li}}
\newcommand{\Lila}{{\rm Le}}
\newcommand{\coez}{\kappa_0}
\newcommand{\coe}{\kappa}
\newcommand{\coeprz}{\kappa'_0}
\newcommand{\coepr}{\kappa'}
\newcommand{\Imagi}{{\rm Im}}
\newcommand{\Reelle}{{\rm Re}}

\newcommand{\Dem}{\noindent{\sc Preuve} : }
\newcommand{\findem}{\end{pf}}

\newcommand{\combi}{\binom}
\newcommand{\equivalent}{\sim}

\newcommand{\gdunrec}{U_{{\rm R}, n}}
\newcommand{\gdvnrec}{V_{{\rm R}, n}}
\newcommand{\ukrec}{u_{{\rm R}, k}}
\newcommand{\unrec}{u_{{\rm R}, n}}
\newcommand{\vnrec}{v_{{\rm R}, n}}
\newcommand{\unmurec}{u_{{\rm R}, n-1}}
\newcommand{\ukmurec}{u_{{\rm R}, k-1}}
\newcommand{\vnmurec}{v_{{\rm R}, n-1}}
\newcommand{\uzerorec}{u_{{\rm R}, 0}}
\newcommand{\uunrec}{u_{{\rm R}, 1}}
\newcommand{\vzerorec}{v_{{\rm R}, 0}}
\newcommand{\vunrec}{v_{{\rm R}, 1}}

\newcommand{\unexpli}{u_{{\rm E},n}}
\newcommand{\vnexpli}{v_{{\rm E},n}}
\newcommand{\unsoro}{u_{{\rm P},n}}
\newcommand{\vnsoro}{v_{{\rm P},n}}
\newcommand{\unnes}{u_{\Sigma,n}}
\newcommand{\vnnes}{v_{\Sigma,n}}
\newcommand{\unbeu}{u_{\R,n}}
\newcommand{\vnbeu}{v_{\R,n}}
\newcommand{\unball}{u_{{\rm TB},n}}
\newcommand{\vnball}{v_{{\rm TB},n}}
\newcommand{\unmod}{u_{{\rm M},n}}  
\newcommand{\vnmod}{v_{{\rm M},n}}  
\newcommand{\vzeromod}{v_{{\rm M},0}}  
\newcommand{\uin}{u_{i,n}}
\newcommand{\vin}{v_{i,n}}
\newcommand{\Iin}{I_{i,n}}
\newcommand{\Ijn}{I_{j,n}}

\newcommand{\hhh}{{\mathfrak{H}}}

\newcommand{\Ibeu}{I_{\R,n}}
\newcommand{\Ines}{I_{\Sigma,n}}
\newcommand{\Icompl}{I_{\C,n}}
\newcommand{\Iball}{I_{{\rm TB},n}}
\newcommand{\Isoro}{I_{{\rm P},n}}
\newcommand{\Irv}{J_n}
\newcommand{\Iriv}{K_n}

\newcommand{\ensembleindices}{\{{\rm R}, {\rm E}, \R, \Sigma,\C,{\rm P},
					{\rm TB},{\rm M} \}}

\newcommand{\fun}{G_n}
\newcommand{\fde}{F_n}

\newcommand{\zeron}{\{0,\ldots,n\}}
\newcommand{\una}{\{1,\ldots,a\}}
\newcommand{\zeroa}{\{0,\ldots,a\}}

\newcommand{\prfrac}[2]{%
 \vcenter{\tabskip0pt\offinterlineskip\halign{%
 \strut##&\hfill\hskip1pt##\hskip1pt\hfill&##\cr
 &$#1$&\vrule\cr\noalign{\hrule}\vrule&$#2$&\cr}}%
}

\newcommand{\vp}{{\rm v}_p}

\newcommand{\del}{\delta}
\newcommand{\cercle}{C}

\newcommand{\congru}{\equiv}
\newcommand{\al}{\alpha}

\newcommand{\R}{\mathbb{R}}

\newcommand{\Q}{\mathbb{Q}}
\newcommand{\Z}{\mathbb{Z}}
\newcommand{\C}{\mathbb{C}}
\newcommand{\scinq}{{\mathfrak S}_5}
\newcommand{\psd}{\rtimes}

\begin{document}
\maketitle

\noindent{\bf INTRODUCTION}

\medskip

Cet exposé est consacr\'e aux valeurs aux entiers $s \geq 2$
de la fonction z\^eta de Riemann, d\'efinie par $\zeta(s) = \sum_{n=1} ^\infty 
n^{-s}$. Quand $s = 2k$ est pair, on sait que $\zeta(2k) \pi^{-2k}$ est
un nombre rationnel, li\'e aux nombres de Bernoulli.
Comme $\pi$ est transcendant (voir  l'appendice de
\cite{Lang} pour une preuve), 
$\zeta(2k)$ l'est aussi pour tout $k \geq 1$. La nature
arithm\'etique des $\zeta(2k+1)$ est beaucoup moins bien
connue. D'un point de vue
conjectural, la situation est simple~:
\begin{conj} \label{conjun}
Les nombres $\pi$, $\zeta(3)$, $\zeta(5)$,  $\zeta(7)$, \ldots sont 
alg\'ebriquement
ind\'ependants sur $\Q$. 
\end{conj}   

Cette conjecture est un cas particulier d'une conjecture diophantienne
sur les polyz\^etas (voir \cite{MiW} ou \cite{CartierMZV}). 
Elle implique que les $\zeta(2k+1)$ 
sont tous transcendants, donc irrationnels, et lin\'eairement ind\'ependants
sur $\Q$.

Tr\`es peu de r\'esultats sont connus en direction de la conjecture~\ref{conjun}. Le premier d'entre eux a \'et\'e annonc\'e par
Ap\'ery  lors des Journ\'ees Arithm\'etiques de Luminy, en 1978~:
\begin{theo}[\cite{Apery}] $\zeta(3)$ est irrationnel.
\end{theo}
Ap\'ery lui-m\^eme n'a donn\'e lors de son expos\'e
(voir \cite{MendesFrance}), et n'a publi\'e \cite{Apery},
qu'une esquisse de sa preuve.
Les d\'etails (qui sont loin d'\^etre triviaux) ont \'et\'e publi\'es 
par Van   
Der Poorten \cite{VDP} (voir aussi \cite{CohenGrenoble} et
\cite{Reyssat}), grâce \`a des contributions de Cohen 
et Zagier.  Par la suite,
plusieurs autres d\'emonstrations du th\'eor\`eme d'Ap\'ery sont parues. La premi\`ere
partie de ce texte est consacr\'ee \`a une synth\`ese des diff\'erents points de
vue qu'on peut adopter pour le d\'emontrer.

La grande perc\'ee suivante date de 2000~:  
\begin{theo}[\cite{RivoalCRAS}, \cite{BR}] \label{thinfinite}
Le $\Q$-espace vectoriel
engendr\'e par $1$, $\zeta(3)$, $\zeta(5)$, $\zeta(7)$, \ldots est de dimension
infinie.
\end{theo}
En cons\'equence, il existe une infinit\'e de $k$ tels que $\zeta(2k+1)$ soit
irrationnel.
On peut donner
des versions effectives de ce dernier \'enonc\'e~: Rivoal a d\'emontr\'e
\cite{vingtetun} que parmi les neuf nombres $\zeta(5)$, 
$\zeta(7)$,  \ldots,  $\zeta(21)$,  l'un au moins est irrationnel.
Ce r\'esultat a \'et\'e am\'elior\'e par Zudilin~:
\begin{theo}[\cite{Zudilinonze}, \cite{Zudilincinqaout}] \label{theoonze}
L'un au moins des quatre nombres $\zeta(5)$,
$\zeta(7)$,  $\zeta(9)$,  $\zeta(11)$ est irrationnel.
\end{theo}
Malgr\'e ces d\'eveloppements r\'ecents, il n'existe aucun entier $s \geq 5$ impair
pour lequel on sache  si $\zeta(s)$ est rationnel ou non. 

\smallskip

Ce texte est divis\'e en trois parties. La premi\`ere est une synth\`ese des
m\'ethodes connues pour d\'emontrer l'irrationalit\'e de $\zeta(3)$~;
l'int\'er\^et des diff\'erentes approches est qu'elles se g\'en\'eralisent plus
ou moins facilement \`a d'autres situations.
La
deuxi\`eme partie fournit une preuve du th\'eor\`eme \ref{thinfinite}, et de
r\'esultats voisins. La troisi\`eme est consacr\'ee \`a des r\'esultats
``quantitatifs''~: mesure d'irrationalit\'e de $\zeta(3)$ et th\'eor\`eme
\ref{theoonze}.

\medskip

\textsc{Remerciements :} Je remercie toutes les personnes qui m'ont aid\'e
dans la pr\'eparation de ce texte, notamment
F.~Amoroso,
V.~Bosser,
N.~Brisebarre,
P.~Cartier, 
G.~Christol,
P.~Colmez,
P.~Grinspan,
L.~Habsieger,
M.~Huttner,
C.~Krattenthaler,
C.~Maclean,
F.~Martin,
Yu.~Nesterenko,
F.~Pellarin,
A.~Pulita,
E.~Royer,  
M.~Waldschmidt, 
D.~Zagier
et W.~Zudilin. Je remercie tout particuli\`erement 
T. Rivoal pour les nombreuses
discussions tr\`es instructives que nous avons eues.


\mathversion{bold}
\section{Irrationalit\'e de $\zeta(3)$}
\mathversion{normal}




Toutes les preuves connues de l'irrationalit\'e de $\zeta(3)$ ont la 
m\^eme structure. 
On 
construit,
pour tout $n \geq 0$, des nombres
rationnels $u_n$ et $v_n$ 
ayant les propri\'et\'es suivantes~:
\begin{enumerate}
\item \label{pointun} La forme lin\'eaire $I_n = u_n \zeta(3) - v_n$ v\'erifie 
$$\limsup_{n \to \infty} \vert I_n \vert^{1/n} \leq
(\sqrt{2}-1)^4 = 0,0294372\ldots$$
\item \label{pointdeux} En notant $d_n$ le p.p.c.m. des entiers compris entre 1 et $n$, 
les coefficients  $u_n$ et $v_n$  v\'erifient~:
$$u_n \in \Z \mbox{ et } 2 d_n ^3 v_n \in \Z.$$
\item\label{pointtrois} Pour une infinit\'e d'entiers $n$, on a $I_n \neq 0$.
\end{enumerate}
La conclusion est alors imm\'ediate~: si $\zeta(3)$ \'etait un nombre
rationnel $p/q$,
alors $2q d_n ^3 I_n$ serait un entier pour tout $n$, et tendrait vers z\'ero
quand $n$ tend vers l'infini (car
$(\sqrt{2}-1)^4 e^3 < 1$, en utilisant 
\cite{Ingham} le
th\'eor\`eme des nombres premiers sous la forme $\lim_{n \to
\infty} \frac{\log(d_n)}{n} = 1$)~: cela contredit la troisi\`eme assertion.

\begin{rema} Comme $(\sqrt{2}-1)^4 \cdot 3,23^3 < 1$,
le th\'eor\`eme des nombres premiers peut \^etre
remplac\'e par l'assertion plus faible $d_n < 3,23 ^n$ pour $n$ assez
grand, qui
se d\'emontre en utilisant des arguments \'el\'ementaires 
\`a la Tchebychev (\cite{Niven}, \S 8.1~; \cite{Ingham}, p. 15). 
\end{rema}


Dans la suite, on donne plusieurs 
constructions (\S~\ref{subsecrec} \`a~\ref{subsecmod})
de $u_n$, $v_n$ et $I_n$, \`a chaque fois 
not\'ees $\uin$, $\vin$ et $\Iin$ (l'indice 
$i \in \ensembleindices$ fait r\'ef\'erence \`a la
construction utilis\'ee). En fait, on construit 
toujours les
m\^emes formes lin\'eaires~: {\em a posteriori} on s'aper\c{c}oit que 
$\uin$, $\vin$ et $\Iin$ ne d\'ependent pas de $i$. La preuve
de cette ind\'ependance est le plus souvent directe. Parfois, on montre
simplement que $\Iin = \Ijn$~; les deux autres \'egalit\'es en
d\'ecoulent en utilisant l'irrationalit\'e de $\zeta(3)$. 

Les premi\`eres valeurs de $u_n$ et $v_n$ sont~: 
\begin{eqnarray*}
(u_n)_{n \geq 0} &=& 1, 5, 73, 1445, 33001, 819005, \ldots \\
(v_n)_{n \geq 0} &=& 0, 6, \frac{351}{4}, \frac{62531}{36}, 
				\frac{11424695}{288}, \ldots 
\end{eqnarray*}

\smallskip

Cette partie contient l'esquisse de plusieurs preuves de l'irrationalit\'e de
$\zeta(3)$, notamment celles d'Ap\'ery \cite{Apery} (\S \ref{subsecrec} et 
\ref{subsecexpli}), de Beukers \cite{Beukers} 
par les int\'egrales multiples (\S \ref{subsecbeukers}) ou \cite{BeukersBesancon}
par les formes modulaires (\S \ref{subsecmod}), de Prevost \cite{Prevost}
(\S \ref{subsecrec} et 
\ref{subsecexpli}),
de Nesterenko \cite{Nesterenko96} (\S \ref{subsecnesterenko}
et \ref{subseccplx}), de Sorokin \cite{SorokinApery} (\S \ref{subsecsorokin}),
et de nombreuses variantes. Certaines preuves sont obtenues en montrant que
deux constructions diff\'erentes fournissent les m\^emes formes lin\'eaires,
puis en prouvant le point \eqref{pointdeux} \`a l'aide de l'une et 
les points \eqref{pointun} et \eqref{pointtrois} \`a l'aide de l'autre (par 
exemple en montrant que $\lim_{n \to \infty} \vert I_n \vert ^{1/n} =
(\sqrt{2}-1)^4$).

\bigskip

La plupart des m\'ethodes connues pour d\'emontrer 
des r\'esultats d'irrationalit\'e
sur les valeurs de $\zeta$  sont li\'ees aux 
polylogarithmes, d\'efinis
pour tout entier $k \geq 1$ par~:
$$\Li_k (z ) = \sum_{n=1} ^\infty \frac{z^n}{n^k},$$
avec $\vert z \vert < 1$ si $k = 1$ et $\vert z \vert \leq 1$ si $k \geq 2$. 
L'id\'ee est de construire des formes lin\'eaires en polylogarithmes, \`a 
coefficients polynomiaux, puis de
sp\'ecialiser en $z=1$. C'est la m\'ethode employ\'ee
dans les  paragraphes~\ref{subsecbeukers} \`a~\ref{subsecball}. 
Les formes lin\'eaires en polylogarithmes $\Iin (z)$ qu'on utilise
ne sont pas toujours les m\^emes, mais elles co\"{\i}ncident en
$z=1$, pour donner les formes lin\'eaires d'Ap\'ery.

Les polylogarithmes s'ins\`erent dans la famille des s\'eries
hyperg\'eom\'etriques $_{q+1}F_q$  (avec $q \geq 1$), d\'efinies par~:
\begin{eqnarray*}
{}_{q+1}F_q
\left(
\begin{array}{cccc}
\alpha_0,&\alpha_1,&\ldots,&\alpha_{q}\\
& \beta_1,&\ldots,&\beta_q\\
\end{array}
\bigg\arrowvert z \right)=\sum_{k=0}^{\infty}
\frac{(\alpha_0)_k(\alpha_1)_k\cdots(\alpha_{q})_k}
{k! \, (\beta_1)_k\cdots(\beta_q)_k} z^k\;,
\label{eq:hyper}
\end{eqnarray*}
o\`u le symbole de Pochhammer est  
$(\al)_k=\alpha(\alpha+1)\cdots(\alpha+k-1)$. 
Dans cet expos\'e, les $\alpha_j$ et les $\beta_j$ seront des entiers,
les $\beta_j$ \'etant positifs,
et $z$ sera un nombre complexe avec $\vert z \vert \leq 1$.
On adopte les d\'efinitions suivantes 
(\cite{AAR}, \S 3.3 et 3.4)~:
\begin{itemize}
\item ${}_{q+1}F_q$ est dite \textit{bien \'equilibr\'ee} si
$\alpha_0+1=\al_1+\beta_1=\cdots=\alpha_{q}+\beta_q$~;
\item ${}_{q+1}F_q$ est dite \textit{tr\`es bien \'equilibr\'ee} si elle est 
bien \'equilibr\'ee et $\alpha_1=\frac12\alpha_0+1$.
\end{itemize}

\subsection{R\'ecurrence lin\'eaire} \label{subsecrec}

\begin{defi} \label{defirec}
Soient $(\unrec)_{n \geq 0} $ et $(\vnrec)_{n \geq 0}$ les 
suites d\'efinies par la relation de r\'ecurrence
\begin{equation} \label{relnrec}
(n+1)^3 y_{n+1} - (34 n^3 + 51 n^2 + 27 n +5) y_n + n^3 y_{n-1} = 0
\end{equation} 
et les conditions initiales
$$\uzerorec = 1 \mbox{ , } \uunrec = 5  \mbox{ , } \vzerorec = 0 \mbox{ , } 
\vunrec  = 6.$$
\end{defi}
Une r\'ecurrence imm\'ediate montre que les suites 
$(\unrec)$ et $(\vnrec)$ sont croissantes et \`a termes rationnels,
avec $n!^3 \unrec \in \Z$ et $n!^3 \vnrec \in \Z$. 
En fait on verra qu'on peut remplacer $n!^3$ par $d_n ^3$.

\medskip

Les propri\'et\'es asymptotiques des suites v\'erifiant la r\'ecurrence
\eqref{relnrec} sont faciles \`a d\'eterminer (voir par exemple 
\cite{Gelfond}, Chapitre 5). L'\'equation caract\'eristique associ\'ee 
est
$X^2 - 34 X + 1$~; elle a deux racines simples, $(\sqrt{2}+1)^4$ et 
$(\sqrt{2}-1)^4$. L'espace vectoriel des solutions de \eqref{relnrec}
est de dimension deux, et admet une base form\'ee de suites
$(y_n ^{(0)})_{n \geq 0}$
 et $(y_n ^{(1)})_{n \geq 0}$ avec 
$\lim_{n \rightarrow + \infty} \frac{\log \vert y_n ^{(0)} \vert}{n} 
 = \log((\sqrt{2}+1)^4)$
 et 
$\lim_{n \rightarrow + \infty} \frac{\log \vert y_n ^{(1)} \vert}{n}  
 = \log((\sqrt{2}-1)^4)$.
La suite $(y_n ^{(1)})$ est uniquement d\'etermin\'ee (\`a proportionnalit\'e
pr\`es) par son comportement asymptotique  ~; toutes les autres
solutions de \eqref{relnrec} se comportent comme $(y_n ^{(0)})$. 
Comme $(\unrec)$ et $(\vnrec)$ sont croissantes, on a :
\begin{equation} \label{asyunvnrec}
\lim_{n \to \infty} \unrec ^{1/n} =
\lim_{n \to \infty} \vnrec ^{1/n} =
(\sqrt{2}+1)^4 = 33,9705627\ldots
\end{equation} 
Quand on adopte ce point de vue, on a int\'er\^et \cite{VDP}
\`a consid\'erer
$\Delta_n =  \begin{array}{|cc|} \vnrec & \vnmurec \\ \unrec & \unmurec 
\end{array}$ pour $n \geq 1$.
La relation de r\'ecurrence montre qu'on a $\Delta_n  = \frac{6}{n^3}$ pour 
tout $n$, ce qui signifie 
$\frac{\vnrec}{\unrec} - \frac{\vnmurec}{\unmurec} = \frac{6}{n^3 \unrec
\unmurec}$.
Donc la suite $(\frac{\vnrec}{\unrec})$ est strictement croissante
et tend vers une limite finie $\ell$, avec 
$ \unrec \ell - \vnrec = \sum_{k = n+1} ^\infty \frac{6 \unrec}{k^3
\ukrec \ukmurec}$.
Ceci prouve que 
$\unrec \ell - \vnrec$ est une solution de \eqref{relnrec} qui tend vers z\'ero
quand $n$ tend vers l'infini~: son comportement asymptotique est
n\'ecessairement donn\'e par
$$\lim_{n \rightarrow +\infty} \frac{\log \vert \unrec \ell - \vnrec \vert}{n}
 = \log((\sqrt{2} - 1)^4).$$

\smallskip

Avec cette d\'efinition de $\unrec$ et $\vnrec$, il n'est pas \'evident
de d\'emontrer que $\ell = \zeta(3)$, et de borner par $d_n^3$
les d\'enominateurs de $\unrec$ et $\vnrec$. Pour ceci, une possibilit\'e 
est de faire le lien avec le paragraphe~\ref{subsecexpli}~: c'est la m\'ethode
employ\'ee dans les premi\`eres preuves d\'etaill\'ees de l'irrationalit\'e de
$\zeta(3)$, qui sont parues peu apr\`es l'expos\'e d'Ap\'ery
(\cite{Reyssat}, \cite{VDP}, \cite{CohenGrenoble}).
 
\medskip


\begin{rema}
Le raisonnement ci-dessus montre que $\frac{\vnrec}{\unrec}$ est la 
$n$-i\`eme  somme
partielle de la s\'erie $\zeta(3) = \sum_{k= 1}^{\infty} \frac{6}{k^3 
\ukrec \ukmurec}$. 
\end{rema}

La d\'efinition~\ref{defirec} s'interpr\`ete en termes de fractions continues
g\'en\'eralis\'ees. En effet, consid\'erons la r\'ecurrence lin\'eaire
\begin{equation} \label{relnrecgdy}
Y_{n+1} - (34 n^3 + 51 n^2 + 27 n +5) Y_n + n^6 Y_{n-1} = 0.
\end{equation} 
On passe d'une solution de \eqref{relnrec} \`a une solution de 
\eqref{relnrecgdy}, et r\'eciproquement, en posant 
$Y_n = n!^3 y_n$. Si
$\gdunrec$ et $\gdvnrec$ sont ainsi associ\'ees \`a $\unrec$ et $\vnrec$, alors
$\frac{\gdvnrec}{\gdunrec} = \frac{\vnrec}{\unrec}$ est la $n$-i\`eme 
r\'eduite de la fraction continue g\'en\'eralis\'ee 
$$\zeta(3)=\prfrac{6 \,}{\, 5}-\prfrac1{117}-\prfrac{64}{\, 535}-\cdots
-\prfrac{n^6}{\, 34n^3+51n^2+27n+5}-\cdots.$$
On peut trouver cette formule gr\^ace \`a 
un proc\'ed\'e g\'en\'eral (\cite{AperyBNF}, \cite{BatutOlivier},
\cite{Zeilbergerdeconstruction}) qui
acc\'el\`ere la convergence d'un d\'eveloppement en
fraction continue g\'en\'eralis\'ee. Ce proc\'ed\'e  
s'applique, en particulier, au d\'eveloppement dont les r\'eduites sont les
sommes partielles de la s\'erie
$\sum_{n=1} ^\infty
\frac{1}{f(n)}$, o\`u $f$ est un polyn\^ome sans z\'ero parmi les
entiers strictement
positifs. 

En utilisant cette 
m\'ethode d'acc\'el\'eration de convergence,  
Andr\'e-Jeannin a d\'emontr\'e \cite{AndreJeannin} que la somme 
des inverses des
nombres de Fibonacci est irrationnelle (voir aussi \cite{BundschuhVaananen}
et \cite{PrevostFibonacci}).

\subsection{Formules explicites} \label{subsecexpli}

\begin{defi} Soient $(\unexpli) $ et $(\vnexpli)$ les 
suites d\'efinies par les formules suivantes~:
\begin{eqnarray*}
\unexpli &=& \sum_{k=0} ^n \combi{n}{k}^2 \combi{n+k}{k}^2 \\
\vnexpli &=& \sum_{k=0} ^n  \combi{n}{k}^2 \combi{n+k}{k}^2
\left( \sum_{m=1} ^n \frac{1}{m^3} + \sum_{m=1} ^k 
\frac{(-1)^{m-1}}{2 m^3 \combi{n}{m} \combi{n+m}{m}} \right)
\end{eqnarray*}
\end{defi}
Sous cette forme, il est clair que $\unexpli \in \Z$ et que 
$\frac{\vnexpli}{\unexpli}$ tend vers $\zeta(3)$. Pour d\'emontrer
(\cite{VDP}, \cite{CohenGrenoble}, \cite{Reyssat})
que $2 d_n ^3 \vnexpli \in \Z$, il suffit de d\'emontrer que,
pour $1 \leq m \leq k \leq n$, 
\begin{equation} \label{quotiententier}
\frac{\combi{n+k}{k} d_n ^3}{m^3 \combi{n}{m} \combi{n+m}{m}}
 = \frac{\combi{n+k}{k-m} d_n ^3}{m^3 \combi{n}{m} \combi{k}{m}}
 \end{equation}
est entier. Soit $p$ un nombre premier~; la valuation $p$-adique $\vp(n!)$
de $n!$ vaut $\sum_{i=1} ^\al
[\frac{n}{p^i}]$ avec $\al = [ \frac{\log(n)}{\log(p)}]
 = \vp(d_n)$. Pour $1 \leq i \leq \vp(m)$ on a 
 $[\frac{n}{p^i}] = [\frac{n-m}{p^i}] + [\frac{m}{p^i}] $
 et pour $\vp(m) < i \leq \vp(d_n)$ on a 
  $[\frac{n}{p^i}] \leq [\frac{n-m}{p^i}] + [\frac{m}{p^i}] +1 $.
On en d\'eduit  $\vp (\combi{n}{m}) \leq \vp(d_n) - \vp(m)$ et
$\vp(\combi{k}{m}) \leq \vp(d_k) - \vp(m)$. Il en r\'esulte que
$\frac{d_n ^3}{m^3 \combi{n}{m} \combi{k}{m}}$
est un entier, et le quotient \eqref{quotiententier} aussi. 


\bigskip

Montrons maintenant (\cite{VDP}, \cite{CohenGrenoble})
que les suites $(\unexpli)$ et 
$(\vnexpli)$ v\'erifient la
r\'ecurrence \eqref{relnrec}. On pose 
$\lam_{n,k}= \combi{n}{k}^2 \combi{n+k}{k}^2$ pour $k,n \in \Z$, et
$$\BBB_{n,k}=4(2n+1)(k(2k+1)-(2n+1)^2) \lam_{n,k},$$ 
avec les conventions habituelles (i.e. $\lam_{n,k}=0$ si $k < 0$ ou $k >n$).
On a alors
$$\BBB_{n,k} - \BBB_{n,k-1} = 
(n+1)^3 \lam_{n+1,k}- (34 n^3 + 51 n^2 + 27 n +5) 
 \lam_{n,k}  + n^3 \lam_{n-1,k}. $$
En sommant sur $k$, on obtient que la suite $(\unexpli)$ satisfait \`a la
r\'ecurrence \eqref{relnrec}. Pour la suite $(\vnexpli)$, on peut
faire de m\^eme en utilisant la suite double
$$\AAA_{n,k} = \BBB_{n,k} \left( \sum_{m=1} ^n \frac{1}{m^3} + \sum_{m=1} ^k 
\frac{(-1)^{m-1}}{2 m^3 \combi{n}{m} \combi{n+m}{m}} \right) +
\frac{5 (2n+1) k(-1)^{k-1}}{n(n+1)} \combi{n}{k}  \combi{n+k}{k}.$$ 
Ceci d\'emontre qu'on a $\unexpli = \unrec$ et 
$\vnexpli = \vnrec$  pour tout $n \geq 0$. Compte tenu des r\'esultats
d\'emontr\'es au paragraphe \ref{subsecrec}, on obtient une preuve de
l'irrationalit\'e de $\zeta(3)$.

\smallskip

La d\'emonstration donn\'ee ci-dessus que $(\unexpli)$ et 
$(\vnexpli)$ v\'erifient la r\'ecurrence \eqref{relnrec}
n'est qu'une simple v\'erification, \`a condition d'\^etre
capable d'exhiber les suites doubles $\BBB_{n,k}$ et $ \AAA_{n,k}$, ce qui 
n'a pas \'et\'e une t\^ache facile (voir \cite{VDP}, \S 7).
Motiv\'es par ce probl\`eme, 
plusieurs auteurs (notamment Zeilberger) ont ensuite mis au point 
des algorithmes permettant d'exhiber de telles suites doubles. On a ainsi
un moyen automatique de produire des
preuves d'identit\'es (voir \cite{CartierZeilberger}, \cite{Zeilberger}, 
\cite{AegaleB}). De plus, ces preuves sont 
imm\'ediatement v\'erifiables \`a la main. 

\bigskip

Dans les formules ci-dessus, un r\^ole central est jou\'e 
par la suite double  
$c_{n,k} =   \sum_{m=1} ^n \frac{1}{m^3} + \sum_{m=1} ^k 
\frac{(-1)^{m-1}}{2 m^3 \combi{n}{m} \combi{n+m}{m}}$ 
(d\'efinie pour $ 0 \leq k \leq n$).
Elle tend vers $\zeta(3)$ quand $n$ tend vers l'infini, uniform\'ement en $k$.
On a 
$c_{n,n} - c_{n-1,n-1} = \frac{5}{2}\frac{(-1)^{n-1}}{n^3 \combi{2n}{n}}$
et $\lim_{n \to \infty} c_{n,n} = \zeta(3)$ donc~:
\begin{equation} \label{atrois}
\zeta(3) = \frac{5}{2} \sum_{n=1} ^\infty \frac{(-1)^{n-1}}{n^3 
\combi{2n}{n}}.
\end{equation}
Cette s\'erie n'est pas utilis\'ee dans la preuve de l'irrationalit\'e de
$\zeta(3)$, mais elle a un int\'er\^et 
non n\'egligeable puisque les $c_{n,k}$ sont
au c\oe ur des formules explicites d\'efinissant $\unexpli$ et $\vnexpli$. 
C'est pourquoi plusieurs auteurs  
ont cherch\'e des g\'en\'eralisations de
\eqref{atrois} (voir par exemple \cite{VDP}, \cite{VDPQueens}, \cite{CohenSMF},
\cite{Koecher}, \cite{Leshchiner}, \cite{BorweinBradley},
\cite{AlmkvistGranville}), parmi lesquelles 
$\zeta(5) = \frac{5}{2} \sum_{n \geq 1} \frac{(-1)^n}{n^3 \binom{2n}{n}}
\left( \sum_{j=1} ^{n-1} \frac{1}{j^2} - \frac{4}{5n^2} \right)$.
Mais aucune de ces g\'en\'eralisations n'a permis d'obtenir de nouveau
r\'esultat d'irrationalit\'e~: la croissance des 
d\'enominateurs est  trop rapide par rapport \`a la convergence.

\bigskip

Pr\'evost a montr\'e \cite{Prevost} comment interpr\'eter les formules
explicites donn\'ees dans ce paragraphe en termes d'approximants de Pad\'e.
Posons  $\pre(\yg) = \sum_{k \geq 1} \frac{1}{(k+\yg)^3}$, c'est-\`a-dire
$\zeta(3, 1+\yg)$ o\`u $\zeta$ est la fonction z\^eta d'Hurwitz (voir
\cite{WhittakerWatson}, Chapitre XIII). Pour tout $n \geq 1$, consid\'erons les
polyn\^omes suivants :
\begin{eqnarray*}
\Ppre_{n}(\yg) &=& \sum_{k=0} ^n \combi{n}{k}\combi{n+k}{k}
\combi{\yg}{k}\combi{\yg+k}{k} = \, \, \fquatretroispre{1} \\
\mbox{ et } \, \, \, \, 
\Qpre_{n}(\yg) &=& \sum_{k=0} ^n \combi{n}{k}\combi{n+k}{k}
\combi{\yg}{k}\combi{\yg+k}{k} \sum_{m=1} ^k 
\frac{(-1)^{m-1}}{2 m^3 \combi{\yg}{m} \combi{\yg+m}{m}}.
\end{eqnarray*}
Alors $\Ppre_{n}$ est de degr\'e $2n$, $\Qpre_{n}$ de degr\'e $2n-2$, 
et on a 
$\Ppre_{n}(\yg) \pre(\yg) - \Qpre_{n}(\yg) = \gdo (\yg^{-2n-1})$
quand $\yg$ tend vers l'infini. Cela signifie que $\Ppre_{n}$ et 
$\Qpre_{n}$ sont des approximants de Pad\'e de la fonction $\pre$. 
Quand $\yg$ est un entier $n$, on a $\pre(n) = 
\zeta(3) - \sum_{m=1} ^n \frac{1}{m^3}$ d'o\`u
$\Ppre_{n}(n) \pre(n) - \Qpre_{n}(n) = \unexpli \zeta(3) - \vnexpli$. 
On peut en d\'eduire \cite{Prevost} la majoration  
$\vert \unexpli \zeta(3) - \vnexpli \vert \leq \frac{4 \pi^2}{(2n+1)^2
\unexpli}$. Pour conclure, on a besoin d'une minoration asymptotique de
$\unexpli$ comme celle de la formule~\eqref{asyunvnrec}. Il suffit donc de v\'erifier
que $\unexpli$ satisfait \`a la r\'ecurrence \eqref{relnrec}. 
On peut utiliser $\BBB_{n,k}$ comme
ci-dessous~; une autre m\'ethode \cite{AskeyWilson} est d'utiliser des relations
de contigu\"{\i}t\'e entre s\'eries hyperg\'eom\'etriques balanc\'ees.

En effet, $\unexpli$ s'\'ecrit 
$\fquatretrois{1}$. Une s\'erie hyperg\'eom\'etrique
$_4 F_3 \left(
\begin{array}{cccc|}
 \al_0 , & \al_1,  & \al_2, & \al_3 \\
     &  \beta_1 , &  \beta_2 , &  \beta_3 
\end{array}
\, \,  z \right)$
est dite (\cite{Slater}, \S 2.1.1)
\textit{balanc\'ee} (ou \textit{Saalsch\"utzienne})
si
 $1 + \sum_{i=0} ^3 \al_i= \sum_{j=1} ^3 \beta_j$. 
Si on modifie deux des sept param\`etres d'une  s\'erie balanc\'ee,
en ajoutant ou en retranchant 1 \`a chacun des deux, on peut obtenir \`a nouveau
une s\'erie balanc\'ee. Si c'est le cas, on dit que ces deux s\'eries sont 
\textit{contigu\"{e}s}. Il y a $2 \cdot \combi{7}{2} = 42$ s\'eries 
balanc\'ees qui sont contigu\"es \`a une s\'erie balanc\'ee donn\'ee.
Quand $\al_0$ est un entier n\'egatif (ce qui signifie que la s\'erie
hyperg\'eom\'etrique est en fait un polyn\^ome), il existe des relations
lin\'eaires 
entre les valeurs en 1 de ces 42 s\'eries, dont les coefficients 
sont des polyn\^omes en les param\`etres $\al_0$, \ldots, $\beta_3$ 
(voir \cite{AAR}, \S 3.7).
On peut \cite{AskeyWilson} d\'eduire de ces relations de contigu\"{\i}t\'e
que la suite $\unexpli$ v\'erifie la r\'ecurrence \eqref{relnrec}.


\subsection{Int\'egrale triple r\'eelle} \label{subsecbeukers}

Consid\'erons 
l'int\'egrale suivante, qui a \'et\'e introduite 
par Beukers \cite{Beukers} (voir aussi \cite{BeukersBolyai})~:
$$\Ibeu(z) = \int_0 ^1 \int_0 ^1 \int_0 ^1 \frac{u^n (1-u)^n v^n (1-v)^n w^n 
(1-w)^n}{((1-w)z+uvw)^{n+1}} \dd u \, \dd v \, \dd w. $$
Cette int\'egrale converge pour tout $z \in \C \moins ] - \infty, 0]$. 
Voici une
esquisse de preuve de l'irrationalit\'e de $\zeta(3)$ qui utilise
$\Ibeu(1)$. Les d\'etails se trouvent dans \cite{Beukers}.

\smallskip

Comme le maximum de la fonction
$\frac{u (1-u) v (1-v) w (1-w)}{1-w(1-uv)}$ 
sur le cube unit\'e  vaut
$(\sqrt{2}-1)^4$, on a~:
$$\lim_{n \rightarrow +\infty} \frac{\log (\Ibeu(1)) }{n}
 = \log((\sqrt{2} - 1)^4).$$
Par ailleurs, si on int\`egre $n$ fois par parties par rapport \`a $v$,
qu'on  change $w$ en $\frac{1-w}{1-w(1-uv)}$, et enfin
qu'on int\`egre $n$ fois 
par parties par rapport \`a $u$, on obtient :
$$\Ibeu(1) = \int_0 ^1 \int_0 ^1 \int_0 ^1 
\frac{P_n(u) P_n(v)}{1-w(1-uv)}
 \dd u \, \dd v \, \dd w, $$
o\`u $P_n(X) = \frac{1}{n!} (X^n (1-X)^n)^{(n)}$ est le $n$-i\`eme polyn\^ome de
Legendre. En int\'egrant par rapport \`a $w$, il vient 
$\Ibeu (1)= \int_0 ^1 \int_0 ^1  \frac{- \log(uv)}{1-uv} P_n(u) P_n(v)
 \dd u \, \dd v$.
Or pour tous $k, l \in \{0, \ldots, n\}$ on peut \'ecrire
$\int_0 ^1 \int_0 ^1  \frac{- \log(uv)}{1-uv} u^k v^l
 \dd u \, \dd v = 2 a_{k,l} \zeta(3) + b_{k,l}$
avec $a_{k,l} \in \Z$ et $d_n ^3 b_{k,l} \in \Z$. On a donc~:
$$\Ibeu (1) = 2 (\unbeu \zeta(3) - \vnbeu) \mbox{ avec } \unbeu \in \Z 
\mbox{ et } 2 d_n ^3 \vnbeu \in \Z.$$
Cela termine la preuve de l'irrationalit\'e de $\zeta(3)$.

\subsection{S\'erie de type hyperg\'eom\'etrique} \label{subsecnesterenko}
Posons 
\begin{equation} \label{eqdefr}
R_n(X) = \frac{(X-1)^2 \ldots (X-n)^2}{X^2 (X+1)^2 \ldots (X+n)^2}
 = \frac{(X-n)_n ^2}{(X)_{n+1}^2} = \frac{\Gamma(X)^4}{\Gamma(X-n)^2 
 \Gamma(X+n+1)^2},
\end{equation}
o\`u $\Gamma$ est la fonction Gamma d'Euler, qui v\'erifie
$\Gamma(s+1) = s \Gamma(s)$. En outre, pour $\vert z \vert \geq 1$ on pose~:
\begin{equation} \label{eqdefines}
\Ines (z) = - \sum_{k=1} ^\infty R'_n(k) z^{-k}.
\end{equation}
En suivant 
\cite{BeukersLNM}, \cite{Gutnik83} et \cite{Nesterenko96} 
on  d\'eveloppe  la fraction rationnelle $R_n$ en \'el\'ements simples~:
\begin{equation} \label{eqdcpelsplesines}
R_n(X) = \sum_{i=0} ^n \left( \frac{\al_i}{(X+i)^2} + \frac{\beta_i}{X+i}
\right),
\end{equation}
avec $\al_i =\combi{n}{i}^2  \combi{n+i}{i}^2 $ et
$\beta_i = 2 (-1)^i  \combi{n}{i} \combi{n+i}{i} \sum_{j \in \{0,
\ldots,n\}, j \neq i}
\frac{(-1)^j   \combi{n}{j} \combi{n+j}{j}}{j-i}$ 
pour $i \in \zeron$ (ces formules s'obtiennent en remarquant que
$R_n(X) = (\frac{(X-n)_n}{(X)_{n+1}})^2$~; voir la d\'emonstration du 
lemme~\ref{lemmeNikishin} ci-dessous, ou bien \cite{Colmez},
\cite{Habsieger} ou \cite{Zudilinelementary}).
En utilisant  \eqref{eqdcpelsplesines} pour exprimer 
\eqref{eqdefines} il vient~:
\begin{eqnarray}
\Ines (z)&=& 2 \sum_{i=0} ^n \al_i z^i \sum_{k \geq 1} \frac{z^{-(k+i)}}{(k+i)^3}
+ \sum_{i=0} ^n \beta_i z^i \sum_{k \geq 1} \frac{z^{-(k+i)}}{(k+i)^2} 
					\nonumber \\
&=& 2 A_n(z) \Li_3(1/z) + B_n (z) \Li_2(1/z) + C_n(z) \label{dvlptines}
\end{eqnarray}
o\`u les polyn\^omes $A_n$, $B_n$ et $C_n$ sont d\'efinis par~:
\begin{eqnarray*}
A_n(z) &=& \sum_{i=0} ^n \al_i z^i
= \fquatretrois{z}	 
\\B_n(z) &=& \sum_{i=0} ^n \beta_i z^i 
\\C_n(z) &=& - \sum_{t=0} ^{n-1} z^t \sum_{i=t+1} ^n \left( \frac{2 \al_i}{(i-t)^3}
+ \frac{\beta_i}{(i-t)^2} \right) 
\end{eqnarray*}
Il est clair que les polyn\^omes $A_n(z)$, $d_n B_n(z)$ et $d_n^3 C_n(z)$
sont \`a
coefficients entiers. 
On a $B_n(1) = 0$ car $R_n$ n'a pas de r\'esidu \`a l'infini.
En posant $\unnes = A_n(1)$ et $\vnnes = -C_n(1)/2$ il vient~:
\begin{equation} \label{NestFourier}
\Ines(1) = 
2 (\unnes \zeta(3) - \vnnes) \mbox{ avec } \unnes \in \Z \mbox{ et }
2d_n ^3 \vnnes \in \Z.
\end{equation}  

Pour d\'emontrer l'irrationalit\'e de $\zeta(3)$, il ne reste plus qu'\`a estimer
$\Ines(1)$. On peut le faire en transformant $\Ines(1)$ en une int\'egrale
complexe (voir le paragraphe~\ref{subseccplx})~; c'est ainsi que Nesterenko
d\'emontre \cite{Nesterenko96} le th\'eor\`eme d'Ap\'ery.

\bigskip

On peut d\'emontrer, en utilisant \cite{Zudilinelementary} l'algorithme de 
``creative telescoping'' (\cite{AegaleB}, Chapitre 6),
que $\Ines(1)$, $\unnes$ et
$\vnnes$ satisfont \`a la relation de r\'ecurrence \eqref{relnrec}.
Cela d\'emontre en particulier l'identit\'e $\vnnes=\vnexpli$.

\subsection{Int\'egrale complexe} \label{subseccplx}

Soit $c$ un r\'eel, avec $0 < c < n+1$. 
Pour $z \neq 0$, choisissons une d\'etermination de $\arg(z)$
strictement comprise entre $-2\pi$ et $2 \pi$, et consid\'erons 
l'int\'egrale suivante,
le long de la droite verticale $\Reelle(s) = c$ dans $\C$, orient\'ee de bas 
en haut~: 
\begin{eqnarray}
\Icompl (z) &=& \frac{1}{2i\pi} \int_{c-i \infty} ^{c + i \infty}
\left( \frac{\pi}{\sin(\pi s)} \right)^2 R_n(s) z^{-s} \dd s \nonumber \\
&=& \frac{1}{2i\pi} \int_{c-i \infty} ^{c + i \infty}
\frac{\Gamma(n+1-s) ^2 \Gamma(s) ^4}{\Gamma(n+1+s)^2} z^{-s} \dd s,
						\label{intmeijer}
\end{eqnarray}
cette derni\`ere \'egalit\'e provenant directement
de \eqref{eqdefr} et de la
formule classique
$\frac{\pi}{\sin(\pi s)} = \Gamma(s) \Gamma(1-s)$.
La valeur de $\Icompl(z)$ ne d\'epend pas du choix de $c$ d'apr\`es 
le th\'eor\`eme des r\'esidus.
L'int\'egrale \eqref{intmeijer} est un exemple de $G$-fonction de Meijer
(voir \cite{Luke}, \S 5.2)~:
$$\Icompl(z) = 
G ^{4,2} _{4,4} \left(
\begin{array}{cccc|}
 -n, & -n ,& n+1, & n+1 \\
   0, & 0,  & 0, &  0
\end{array}
\, \, z \right).$$

La m\'ethode du col (voir par exemple \cite{Dieudonne}, Chapitre IX)
permet \cite{Nesterenko96}
d'obtenir une estimation asymptotique tr\`es pr\'ecise~:
$$\Icompl (1)= \frac{\pi ^{3/2} 2^{3/4}}{n^{3/2}} (\sqrt{2}-1)^{4n+2}
(1+\gdo(n^{-1})).$$


\medskip

Quand on d\'eplace le contour d'int\'egration vers la droite
pour faire appara\^{\i}tre les p\^oles $n+1$, $n+2$, \ldots, 
le th\'eor\`eme des r\'esidus donne 
(\cite{Gutnik79}, \cite{Gutnik83}), puisque $(\frac{\pi}{\sin(\pi s)})^2=
\frac{1}{(s-k)^2} + \gdo(1)$ quand $s$ tend vers un entier $k$~:
\begin{equation} \label{eqfunfde}
\Icompl(z) =  \Ines (z) +  \log(z) \sum_{k = 1} ^\infty R_n(k) z^{-k}.
\end{equation} 
En particulier pour $z = 1$ on obtient $\Icompl(1) = \Ines(1)$.

Par ailleurs, Nesterenko 
a d\'emontr\'e \cite{NesterenkoCaen} un th\'eor\`eme g\'en\'eral 
qui relie une int\'egrale multiple
r\'eelle \`a une int\'egrale complexe~; dans notre cas particulier, ce
th\'eor\`eme donne $\Icompl(z) = \Ibeu(z)$. 

On peut d\'emontrer \cite{Nesterenko96} que $\Icompl(1)$
v\'erifie la r\'ecurrence \eqref{relnrec} en 
utilisant les relations de contigu\"{\i}t\'e sur les
$G$-fonctions de Meijer. C'est en fait une preuve parall\`ele \`a celle du
paragraphe~\ref{subsecexpli}, o\`u on utilisait la contigu\"{\i}t\'e entre des
$_4 F_3$. En effet (\cite{Luke}, \S 5.8), 
ces $_4 F_3$ 
satisfont aux m\^emes \'equations
diff\'erentielles que les $G$-fonctions de Meijer correspondantes, 
donc aux m\^emes
relations de contigu\"{\i}t\'e.


\subsection{Un probl\`eme d'approximation de Pad\'e} \label{subsecpade}

Consid\'erons \cite{BeukersLNM}
le probl\`eme suivant~: trouver quatre polyn\^omes 
$A_n$, $B_n$, $C_n$ et $D_n$, \`a
coefficients rationnels, de degr\'e au plus $n$, tels que~:
\begin{equation} \label{pbpade} 
\begin{cases}
\fde (z):=A_n(z)\Li_2(1/z)+B_n(z)\Li_1(1/z)+D_n(z)=\gdo(z^{-n-1}) 
\mbox{ quand } z \to \infty \\ 
\fun (z):=2A_n(z)\Li_3(1/z)+B_n(z)\Li_2(1/z)+C_n(z)=\gdo(z^{-n-1})
 \mbox{ quand } z \to \infty \\
B_n(1) = 0
\end{cases}
\end{equation}  

Une solution \`a ce probl\`eme de Pad\'e est donn\'ee par les polyn\^omes
$A_n$, $B_n$ et $C_n$ du paragraphe~\ref{subsecnesterenko} (et un
polyn\^ome $D_n$ convenable). On a alors~:
$$\begin{cases}
\fde (z) = \sum_{k=1} ^\infty R_n(k) z^{-k} 
 = \frac{n!^4}{(2n+1)!^2} z^{-n-1}
{_4 F_3} \left(
\begin{array}{cccc|}
 n+1, & n+1, & n+1, & n+1 \\
     &  2n+2,  &  2n+2, &  1
\end{array}
\, \, \,  z^{-1} \right)\\
\fun (z) = \Ines (z) = - \sum_{k=1} ^\infty R'_n(k) z^{-k} 
\end{cases} $$
En effet, la seconde \'egalit\'e est simplement une r\'e\'ecriture de 
\eqref{eqdefines} et 
\eqref{dvlptines}. La premi\`ere se d\'emontre de mani\`ere analogue \`a
  \eqref{dvlptines}, mais sans d\'eriver \eqref{eqdcpelsplesines}. 
  
\bigskip

L'\'equation diff\'erentielle hyperg\'eom\'etrique sous-jacente 
aux constructions des paragraphes~\ref{subsecnesterenko} et~\ref{subseccplx}
s'\'ecrit $\ope y = 0$, en posant
$$\ope = z(\de +n +1)^2 (\de -n)^2 - \de^4 \mbox{ avec } 
\de = z \frac{\dd}{\dd z}.$$
Elle admet au voisinage de l'infini quatre solutions lin\'eairement
ind\'ependantes~: $\fde(z)$, $\Icompl(z) = \fun(z) + \fde(z) \log(z)$,
$A_n(z)$ et $B_n(z) - A_n(z) \log(z)$ 
(voir \cite{Luke}, \S 5.1 et
5.8, \cite{Huttnerlille} et \cite{Gutnik83}). 
Ces solutions sont reli\'ees par la monodromie~: en prolongeant
analytiquement $\fde$ le long d'un lacet qui entoure le point 1
on fait appara\^{\i}tre  $B_n(z) + A_n(z) \log(1/z)$, puis  
en faisant le tour de l'infini on obtient $A_n(z)$ (voir 
\cite{Oesterle} pour la monodromie des polylogarithmes). 


Ce point de vue permet de d\'emontrer \cite{Huttnerlille} que le 
probl\`eme de Pad\'e \eqref{pbpade} a une solution unique (\`a 
proportionnalit\'e pr\`es). En effet,
en partant d'une solution 
$A_n$, $B_n$, $C_n$, $D_n$, 
on montre que $\fde$  v\'erifie une \'equation
diff\'erentielle lin\'eaire fuchsienne d'ordre 4
qu'on d\'etermine explicitement (en calculant ses exposants, et en utilisant la
relation de Fuchs)~: on trouve que c'est $\ope y = 0$.

\bigskip

Pour d\'emontrer l'unicit\'e de la solution de ce probl\`eme de Pad\'e,
on peut aussi suivre  \cite{BeukersLNM}.
On part d'une solution quelconque, avec des polyn\^omes  
$A_n$, $B_n$, $C_n$, $D_n$ et des fonctions $\fde$ et $\fun$. On  
note $\al_i$ et $\beta_i$ les coefficients de $A_n$ et $B_n$, 
et on leur associe la fraction
rationnelle $R_n$ d\'efinie  par \eqref{eqdcpelsplesines}. On voit alors que
$\fde (z) = \sum_{k=1} ^\infty R_n(k) z^{-k} $ et
$\fun (z) = - \sum_{k=1} ^\infty R'_n(k) z^{-k}$ ,  donc 
les deux premi\`eres contraintes de \eqref{pbpade} signifient 
que $R_n$ et sa d\'eriv\'ee s'annulent aux points 1, 2, \ldots, $n$. En outre,
le r\'esidu \`a l'infini  de $R_n$ est alors $B_n(1)=0$~: la fraction
rationnelle $R_n$ est n\'ecessairement donn\'ee, \`a constante multiplicative
pr\`es, par \eqref{eqdefr}.

\subsection{Polyn\^omes orthogonaux}

Consid\'erons (\cite{BorweinErdelyi}, \cite{AsscheDelhi})
le probl\`eme suivant : trouver deux polyn\^omes 
$\Aorth_n$ et $\Borth_n$, de degr\'e au plus $n$, tels que :
\begin{equation} \label{pbortho}
\begin{cases}
\int_0 ^1 \left( \Borth_n(x) - \Aorth_n(x) \log(x) \right) x^k \dd x = 0
\mbox{ pour tout } k \in \{ 0, \ldots, n-1\} \\
\int_0 ^1 \left( \Borth_n(x)  - \Aorth_n(x)  \log(x) \right) x^k \log(x) \dd x = 0
\mbox{ pour tout } k \in \{ 0, \ldots, n-1\} \\
\Borth_n(1) = 0
\end{cases}
\end{equation}  
Une solution \`a ce probl\`eme est donn\'ee par les polyn\^omes 
$\Aorth_n$ et $\Borth_n$ d\'efinis par~:
\begin{equation} \label{bnanleg}
\Borth_n(x) - \Aorth_n(x) \log(x) = \int_x ^1 P_n(\frac{x}{t}) 
P_n(t) \frac{\dd t}{t},
\end{equation}
o\`u $P_n$ est le $n$-i\`eme polyn\^ome de Legendre (comme au 
paragraphe~\ref{subsecbeukers}). En effet, on a alors
$\int_0 ^1 \left( \Borth_n(x) - \Aorth_n(x) \log(x) \right) x^k \dd x = 
\left( \int_0 ^1 P_n(u) u^k \dd u \right) ^2$ en posant $u =\frac{x}{t}$.
La premi\`ere condition de \eqref{pbortho} en d\'ecoule 
imm\'ediatement ; la
deuxi\`eme s'obtient apr\`es d\'erivation par rapport \`a $k$.

\smallskip

Comme on a $\Li_j(1/z) = \frac{(-1)^{j-1}}{(j-1)!} \int_0 ^1 \log^{j-1}(x)
\frac{\dd x}{z-x}$ pour tout entier $j \geq 1$, il vient~:
\begin{equation} \label{eqresteavecc}
2 \Aorth_n(z) \Li_3(1/z) + \Borth_n(z)\Li_2(1/z) = 
- \int_0 ^1  \left(\Borth_n(z) - \Aorth_n(z) \log(x) \right) \frac{\log(x) \, \dd x}{z-x}.
\end{equation}
On d\'efinit un polyn\^ome $\Corth_n(z)$ par~:
$$\Corth_n(z) = 
\int_0 ^1  \frac{\Borth_n(z) -\Borth_n(x)}{z-x} \log(x) \, \dd x -
 \int_0 ^1  \frac{\Aorth_n(z) -\Aorth_n(x)}{z-x} \log^2(x) \, \dd x .$$
Gr\^ace \`a \eqref{bnanleg} on peut obtenir des formules explicites pour
$\Aorth_n$, $\Borth_n$ et $\Corth_n$~; on trouve les m\^emes que
pour $A_n$, $B_n$ et $C_n$ respectivement au paragraphe~\ref{subsecnesterenko}.
Donc $\Aorth_n$, $d_n \Borth_n$ et $d_n^3 \Corth_n$ sont \`a
coefficients entiers. On obtient aussi 
(\cite{BorweinErdelyi}, Corollaire A.2.3) que tous les z\'eros
de $\Aorth_n(z)$ et de $\frac{\Borth_n(z)}{z-1}$ sont r\'eels n\'egatifs, et 
entrelac\'es.
Par ailleurs, on a~:
$$2 \Aorth_n(z) \Li_3(1/z) + \Borth_n(z)\Li_2(1/z) + \Corth_n(z) =
- \int_0 ^1  \left(\Borth_n(x) - \Aorth_n(x) \log(x) \right) \frac{\log(x) \dd x}{z-x}.$$
Quand $z=1$, le membre de droite se transforme
(en utilisant \eqref{bnanleg} et en posant $u=t$, 
 $v = \frac{x}{t}$) en
$ \Ibeu(1) = - \int_0 ^1 \int_0 ^1 \frac{\log(uv)}{1-uv}P_n(u)P_n(v) \dd u \dd v$. 
En appliquant l'estimation asymptotique
de $\Ibeu(1)$ obtenue
au paragraphe~\ref{subsecbeukers}, on obtient 
une d\'emonstration de l'irrationalit\'e de
$\zeta(3)$.


\bigskip

En fait un couple $(\Aorth_n, \Borth_n)$ v\'erifie \eqref{pbortho} si, et seulement si, 
il existe $C_n$ et $D_n$ tels que $(\Aorth_n, \Borth_n, C_n, D_n)$ soit une
solution du probl\`eme de Pad\'e \eqref{pbpade}. Plus pr\'ecis\'ement,
la premi\`ere (resp. la deuxi\`eme) assertion de
\eqref{pbpade} \'equivaut \`a la premi\`ere (resp. la deuxi\`eme) assertion de
\eqref{pbortho} (il s'agit d'un fait g\'en\'eral~: voir par exemple
\cite{NikishinSorokin}, Chapitre 4, \S 3.4). 
D\'emontrons-le pour la deuxi\`eme. Soient
$\Gamma$ un chemin qui entoure
le segment $[0,1]$ dans sens direct, et 
$k \in \{ 0, \ldots, n-1\}$. On a :
\begin{small}
$$\frac{1}{2 i \pi} \int_\Gamma z^k \left(
2 \Aorth_n(z) \Li_3(1/z) + \Borth_n(z) \Li_2(1/z) \right) \, \dd z
= - \int_0 ^1 \left( \Borth_n(x) - \Aorth_n(x) \log(x) \right) x^k \log(x) \, \dd x ,$$
\end{small}
d'apr\`es \eqref{eqresteavecc}, 
en intervertissant les deux signes d'int\'egration et
en appliquant le th\'eor\`eme des r\'esidus. 

Il d\'ecoule de ceci que le probl\`eme \eqref{pbortho} a une solution unique
(\`a proportionnalit\'e pr\`es), donn\'ee par $\Aorth_n = A_n$
et $\Borth_n = B_n$.


\subsection{D'autres probl\`emes d'approximation de Pad\'e} 
						\label{subsecsorokin}

Sorokin \cite{SorokinApery} consid\`ere le probl\`eme de Pad\'e suivant~: pour
$n \geq 0$, trouver des polyn\^omes $\Asoro_n$, $\Bsoro_n$, $\Csoro_n$, 
$\Dsoro_n$ de degr\'e au
plus $n$ tels qu'on ait~:
$$\begin{cases}
\Isoro(z):=\Asoro_n(z)\Lila_{2,1}(1/z)+\Bsoro_n(z)\Lila_{1,1}(1/z)+
	\Csoro_n(z) \Li_1(1/z) + \Dsoro_n(z)=\gdo(z^{-n-1}) \\
\hspace{12cm} \mbox{ quand } z \to \infty\\
\Asoro_n(z)\Li_2(1-z)+\Csoro_n(z)=\gdo ((1-z)^{n+1}) \mbox{ quand } z \to 1 \\
\Asoro_n(z)\Li_1(1-z)+\Bsoro_n(z)=\gdo ((1-z)^{n+1}) \mbox{ quand } z \to 1 ,
\end{cases}$$
o\`u pour $s_1, \ldots, s_k \geq 1$ on d\'efinit le polylogarithme multiple
$$\Lila_{s_1,\ldots,s_k}(z) = \sum_{n_1 \geq \ldots \geq n_k \geq 1} 
\frac{z^{n_1}}{n_1^{s_1}\ldots n_k^{s_k}},$$
qui v\'erifie $\Lila_{2,1}(1) = 2 \zeta(3)$ (voir \cite{MiW}).

Sorokin d\'emontre que ce probl\`eme de Pad\'e
admet une solution unique, et qu'\`a
proportionnalit\'e pr\`es elle v\'erifie (pour 
$z \in \C \moins [0,1[$)~:
$$\Isoro(z) = z^{n+1}
\int_0 ^1 \int_0 ^1 \int_0 ^1 \frac{u^n (1-u)^n v^n (1-v)^n w^n 
(1-w)^n}{(z-uv)^{n+1}(z-uvw)^{n+1}} \dd u \, \dd v \, \dd w. $$
Avec cette normalisation, $\Asoro_n$ est 
\`a coefficients entiers (donc aussi $d_n \Bsoro_n$, $d_n ^2 \Csoro_n$ et 
$d_n^3 \Dsoro_n$), d'o\`u~:
$$\Isoro(1) = 2 (\unsoro \zeta(3) - \vnsoro) \mbox{ avec }
\unsoro \in \Z \mbox{ et } 2 d_n ^3 \vnsoro \in \Z.$$
De plus l'expression int\'egrale donne facilement l'estimation asymptotique de 
$\Isoro(1)$~; c'est ainsi que Sorokin d\'emontre l'irrationalit\'e de
$\zeta(3)$. 

\bigskip

Un th\'eor\`eme g\'en\'eral de Zlobin \cite{Zlobin} montre qu'on a 
$$\Isoro(z) = \int_0 ^1 \int_0 ^1 \int_0 ^1 \frac{u^n (1-u)^n v^n (1-v)^n w^n 
(1-w)^n}{(z-w(1-uv))^{n+1}} \dd u \, \dd v \, \dd w, $$
d'o\`u $\Isoro(1) = \Ibeu(1)$. On peut obtenir directement 
ce r\'esultat en
appliquant le changement de variables (\cite{SFCRAS}, \S 2)
d\'efini par $U = 1-w$, $V = \frac{(1-u)v}{1-uv}$
et $W = u$ (et qui v\'erifie $1-W(1-UV) = \frac{(1-u)(1-uvw)}{1-uv}$).

\bigskip

Il existe plusieurs autres
probl\`emes de Pad\'e li\'es \`a $\zeta(3)$~; l'un d'entre
eux \cite{Sorokin94} fait appara\^{\i}tre l'int\'egrale suivante~:
$$\int_0 ^1 \int_0 ^1 \int_0 ^1 \frac{u^n (1-u)^n v^n (1-v)^n w^n 
(1-w)^n}{(z(1-u+uv)-uvw)^{n+1}} \dd u \, \dd v \, \dd w. $$
Le changement de variables qui fixe $u$ et $w$ et change $v$ en 
$\frac{v}{1-u(1-v)}$ transforme cette int\'egrale en 
$$\int_0 ^1 \int_0 ^1 \int_0 ^1 \frac{u^n (1-u)^n v^n (1-v)^n w^n 
(1-w)^n}{(1-uv)^{n+1}(z-uvw)^{n+1}} \dd u \, \dd v \, \dd w. $$

Ces diff\'erents probl\`emes de Pad\'e fournissent tous les formes lin\'eaires
d'Ap\'ery en 1 et $\zeta(3)$, mais ils correspondent \`a des combinaisons
lin\'eaires diff\'erentes de polylogarithmes.

\subsection{S\'erie hyperg\'eom\'etrique tr\`es bien \'equilibr\'ee} 
						\label{subsecball}
On pose~:
\begin{eqnarray*}
\BB_n(X) &=& n!^2 (2X+n) \frac{(X-1)\ldots(X-n)(X+n+1)\ldots(X+2n)}{X^4
 (X+1)^4 \ldots (X+n)^4} \\
 &=& n!^2 (2X+n) \frac{(X-n)_n (X+n+1)_n}{(X)_{n+1}^4}
\end{eqnarray*}
et
$$\Iball (z) = \sum_{k =1} ^\infty \BB_n(k) z^{-k}.$$

La  s\'erie  $\Iball (1)$
a \'et\'e introduite par K. Ball  (voir \cite{survol})
dans le but de r\'epondre \`a une question de
Nesterenko \cite{Nesterenko96}~: trouver une preuve
de l'irrationalit\'e de $\zeta(3)$ analogue \`a celle de
Fourier (\cite{EMS}, Chapitre 2, \S 1.1) 
pour l'irrationalit\'e de $e$. En effet, on peut 
estimer $\Iball (1)$ de mani\`ere
\'el\'ementaire (\cite{Zudilinelementary}, Lemme 4~;
\cite{theseTanguy}, \S 5.1~; 
voir aussi la seconde d\'emonstration du lemme 3 de \cite{BR})~:
$$\lim_{n \to +\infty} \frac{\log(\Iball (1))}{n} = \log((\sqrt{2}-1)^4),$$
ou bien (voir le paragraphe~\ref{subsecdetails})
d\'eduire cette estimation d'une repr\'esentation int\'egrale
de $\Iball(z)$ vue comme  s\'erie 
hyperg\'eom\'etrique tr\`es bien \'equilibr\'ee~:
$$\Iball (z) = z^{-n-1} \frac{n!^7 (3n+2)!}{(2n+1)!^5} \,
\fpqusz{7}{6}{3n+2}{\frac{3}{2}n+2}{n+1}{n+1}{\frac{3}{2}n+1}{2n+2}{2n+2}.$$
De plus, on a $\Iball(z) = P_0(z) + \sum_{j=1} ^4 P_j(z) 
\Li_j(1/z)$ avec des polyn\^omes $P_0$, \ldots, $P_4 \in \Q[z]$ v\'erifiant 
$P_j(z) = (-1)^{j+1} z^4 P_j(1/z)$ pour tout
$j \in \{1,\ldots,4\}$, $P_1(1)=0$ et 
$d_n ^{4-j} P_j(z) \in \Z[z]$ pour tout
$j \in \{0,\ldots,4\}$ (ceci sera g\'en\'eralis\'e au paragraphe~\ref{subsecdetails}). En particulier, on en d\'eduit
\begin{equation} \label{dcpball}
\Iball (1) = 2 (\unball \zeta(3) - \vnball) \mbox{ avec }
2d_n \unball \in \Z \mbox{ et } 2d_n^4 \vnball \in  \Z.
\end{equation}
Mais ceci ne suffit pas \`a d\'emontrer l'irrationalit\'e de $\zeta(3)$,
car $(\sqrt{2}-1)^4 e ^4 > 1$.


\bigskip

Une identit\'e de Bailey (\cite{Zudilincinqaout}, Proposition 2~;
\cite{Slater}, formule (4.7.1.3)) donne $\Iball(1) = \Icompl(1)$. 
Une telle identit\'e ne peut pas avoir lieu pour tout $z$, car
$\Li_4(1/z)$ appara\^{\i}t dans la d\'ecomposition en polylogarithmes de 
$\Iball(z)$ mais pas dans celle de $\Icompl(z)$.
Par ailleurs Zudilin a d\'emontr\'e une identit\'e 
g\'en\'erale (\cite{Zudilinservice}, Th\'eor\`eme 5) qui \'ecrit une 
s\'erie hyperg\'eom\'etrique  tr\`es bien \'equilibr\'ee sous la forme
d'une int\'egrale
g\'en\'eralisant celles introduites par Beukers \cite{Beukers}, 
Vasilenko \cite{Vasilenko} et Vasilyev (\cite{Vasilyevancien}, 
\cite{Vasilyev}). Dans
notre cas particulier, cette identit\'e est $\Iball(1) = \Ibeu(1)$.
Enfin, en utilisant les  algorithmes 
d\'ecrits dans \cite{AegaleB} on peut
d\'emontrer  que 
$\Iball(1)$ (\cite{theseTanguy}, \S 5.1~; \cite{Zudilinelementary}), 
ainsi que $\unball$ et $\vnball$ \cite{Krattenthaler}, 
v\'erifient la relation de
r\'ecurrence \eqref{relnrec}.
On en d\'eduit $\unball = \unexpli$ et $\vnball = \vnexpli$, d'o\`u 
$\unball \in \Z$  et  $2d_n^3 \vnball \in \Z$ (ce qui est plus pr\'ecis
que \eqref{dcpball}).




\subsection{Preuve utilisant des formes modulaires} \label{subsecmod}

Dans ce paragraphe, on esquisse une preuve due \`a Beukers
\cite{BeukersBesancon}
de l'irrationalit\'e de $\zeta(3)$. Les outils mis en \oe uvre sont expos\'es
dans \cite{Serre} (Chapitre VII) et \cite{Zagier}.

Pour $\tau$ dans le demi-plan de Poincar\'e $\hhh$, posons $q = e^{2 i \pi
\tau}$ et consid\'erons les s\'eries d'Eisenstein
$E_2(\tau) = 1 - 24 \sum_{n \geq 1} \sigma_1(n) q^n$ et 
$E_4(\tau) = 1 + 240 \sum_{n \geq 1} \sigma_3(n) q^n$.
On pose~: 
\begin{eqnarray*}
E(\tau) &=& \frac{1}{24} \left( -5 E_2(\tau) + 2 E_2(2 \tau) - 3 E_2 (3\tau)
+ 30 E_2 (6 \tau) \right) \\
\mbox{et \, \, \, }
F(\tau) &=& \frac{1}{40} \left(E_4(\tau) -28  E_4(2 \tau) + 63 E_4 (3\tau)
-36  E_4 (6 \tau) \right).
\end{eqnarray*}
Alors $E(\tau)$, respectivement $F(\tau)$, est une forme modulaire de poids
2, resp. 4, pour $\Gamma_0(6)$. Si $F(\tau) =  \sum_{n \geq 1} f_n q^n$
d\'esigne le d\'eveloppement de Fourier de $F$ \`a l'infini (o\`u elle
s'annule), on pose $f(\tau) =  \sum_{n \geq 1} \frac{f_n}{n^3} q^n$.
On a alors
$(\frac{\dd}{\dd \tau})^3 f(\tau) = (2 i \pi )^3 F(\tau)$.

Consid\'erons la fonction modulaire  pour $\Gamma_0(6)$ donn\'ee par~:
$$t(\tau) \, = \, \left( \frac{\Delta(6\tau) \Delta(\tau)}{\Delta(2\tau)
\Delta(3\tau)} \right) ^{1/2} \, = \, \, \, 
q \! \! \! \! \! \! \prod_{{\tiny \begin{array}{c} n \geq 1 \\ {\rm pgcd}(n,6)=1 \end{array}
}} \! \! \! \! \! \! \! (1-q^n)^{12},$$
avec $\Delta (\tau) = q \prod_{n \geq 1} (1-q^n)^{24}$.
Elle n'a ni z\'ero ni p\^ole dans $\hhh$. Au voisinage de $q=0$,
$t(\tau) = q - 12 q^2 + 66 q^3 - \ldots $ s'\'ecrit 
comme une s\'erie enti\`ere en $q$, \`a coefficients entiers, 
avec un rayon de convergence \'egal \`a 1. 
Elle admet une r\'eciproque locale, not\'ee
$q(t) \in \Z [[t]]$. Par composition, on peut donc d\'efinir des suites 
$(\unmod)$ et $(\vnmod)$ par~:
\begin{eqnarray*}
E(q(t)) &=& \sum_{n \geq 0} \unmod t^n \in \Z [[t]] \\
\mbox{et }E(q(t))f(q(t))  &=& \sum_{n \geq 0} \vnmod t^n \in \Q [[t]]
\mbox{ avec } \vzeromod = 0 \mbox{ et }
d_n ^3 \vnmod \in \Z \mbox{ pour tout } n \geq 1.
\end{eqnarray*}
Notons, pour $k \in \Z$, $w_k$ l'op\'erateur
d'Atkin-Lehner 
d\'efini par 
$(w_k g)(\tau) = 6^{-k/2} \tau^{-k} g(\frac{-1}{6 \tau})$.
Alors $w_2 (E) = -E$ et $w_4(F) = -F$. De cette seconde \'egalit\'e  (et
d'un lemme de Hecke~: voir \cite{Weil}, \S 5) d\'ecoule la relation
$w_{-2}(h) = -h$, en posant $h(\tau) = 
L(F,3) - f(\tau)$, o\`u $L(F,s)$ est la fonction $L$ de $F$. 
Il vient alors $w_0 (Eh) = Eh$, c'est-\`a-dire 
que la fonction 
$E(\tau) h(\tau)$ est invariante par la substitution $\tau \mapsto
\frac{-1}{6 \tau}$. 

Consid\'erons maintenant les rayons de convergence. La fonction $t(\tau)$ est 
ramifi\'ee seulement 
au-dessus des points $(\sqrt{2}-1)^4$, $(\sqrt{2}+1)^4$ et
$\infty$. Au-dessus de $(\sqrt{2}-1)^4$, le seul 
point de ramification
(modulo $\Gamma_0(6)$)
est $\tau = i/\sqrt{6}$~; il est d'indice deux, et les deux branches
en ce point sont
\'echang\'ees par l'involution  $\tau \mapsto
\frac{-1}{6 \tau}$. Comme $E(\tau) h(\tau)$ est invariante par cette involution,
on peut d\'efinir $Eh$ comme une fonction de $t$ au voisinage 
de $ t=  (\sqrt{2}-1)^4$, et en fait sur tout le disque $\vert t \vert 
< (\sqrt{2}+1)^4$. Cela signifie que la s\'erie 
$\sum_{n \geq 0} (L(F,3) \unmod - \vnmod) t^n$ a un rayon de convergence
sup\'erieur ou \'egal \`a $(\sqrt{2}+1)^4$, c'est-\`a-dire qu'on a~:
$$\limsup_{n \to \infty} \frac{\log \vert L(F,3) \unmod - \vnmod \vert}{n}
\leq \log((\sqrt{2}-1)^4).$$
Ceci conclut la d\'emonstration de l'irrationalit\'e de $L(F,3)$. Or on peut
calculer explicitement $L(F,s)$. En effet, quand 
$\Reelle(s) > 4$ on a, pour tout entier $j \geq 1$~:
$$L(E_4(j\tau),s) =1+240 \sum_{n \geq 1} \frac{\sigma_3(n)}{(jn)^s}
 = 1+240 \sum_{d,e \geq 1}  \frac{d^3}{(jde)^s}
= 1+240 \zeta(s) \zeta(s-3) j^{-s}.$$
On en d\'eduit imm\'ediatement $L(F,s) = -2 \zeta(s) \zeta(s-3)$, d'o\`u
$L(F,3) = \zeta(3)$.

\bigskip

Comme $E(\tau)$ est une forme modulaire de poids 2 et $t(\tau)$ une fonction
modulaire, la fonction $E(q(t))$ de la variable $t$ est solution 
\cite{ZagierCDF} (voir aussi \cite{Beukerspideux}, p. 58) 
d'une \'equation diff\'erentielle 
lin\'eaire $\Lmod y=0$, d'ordre 
trois. On peut la d\'eterminer explicitement~:
$$\Lmod = (t^4-34t^3+t^2) \frac{\dd^3}{\dd t^3} + (6t^3-153 t^2 + 3t) 
\frac{\dd^2}{\dd t^2} + (7t^2 - 112 t +1) \frac{\dd }{\dd t} 
+ (t-5).$$
Cette \'equation diff\'erentielle v\'erifi\'ee par la s\'erie g\'en\'eratrice
des $\unmod$ montre qu'ils satisfont \`a la relation de r\'ecurrence
\eqref{relnrec} : on a donc $\unmod= \unrec$ (voir aussi 
\cite{Beukersanother}).
En posant $V(t) = E(q(t))f(q(t))$ on montre \cite{ZagierCDF} que
$\Lmod V = 5$, d'o\`u $\vnmod= \vnrec$.

Une base de solutions de
l'\'equation diff\'erentielle $ \Lmod y=0$ est donn\'ee 
par $E(q(t))$, $\tau(t) E(q(t))$ et $\tau^2(t) E(q(t))$
(voir aussi \cite{BeukersPeters}, Corollaire 2).
La seule solution qui soit r\'eguli\`ere en $0$ est $E(q(t))$ (\`a
proportionnalit\'e pr\`es).
De plus, la construction de $\Lmod$ montre \cite{ZagierCDF} 
que c'est un carr\'e sym\'etrique, ce qui peut se v\'erifier directement 
(voir \cite{DworkAmice}). 



\bigskip

\begin{rema}
Le point de vue adopt\'e dans ce paragraphe est li\'e ``individuellement''
\`a $\zeta(3)$ (qui est vu comme valeur sp\'eciale d'une fonction $L$), 
par opposition aux
m\'ethodes utilis\'ees dans les paragraphes~\ref{subsecbeukers} 
\`a~\ref{subsecball}, o\`u
$\zeta(3)$ apparaissait comme la valeur en 1 d'un polylogarithme. 
\end{rema}


Cette preuve de l'irrationalit\'e de $\zeta(3)$ 
s'exprime naturellement en termes des 
s\'eries g\'en\'eratrices
$U(t) = \sum_{n \geq 0} u_n t^n$ et $V(t) = \sum_{n \geq 0} v_n t^n$
des approximations rationnelles de $\zeta(3)$ (voir \cite{VDPDPP}, 
\cite{BeukersBesancon} et \cite{ChudCarbondale}, \S 5 pour d'autres preuves
dans le m\^eme esprit). 
L'aspect arithm\'etique consiste \`a d\'emontrer que 
les coefficients de $U(t)$ sont entiers, et que
$d_n^3$ est un d\'enominateur commun aux $n$ premiers coefficients de $V(t)$~:
c'est une majoration $p$-adique de ces coefficients, pour toute place 
finie $p$. L'aspect analytique est une minoration, par
$(1+\sqrt{2})^4$, du rayon de convergence (archim\'edien) de la s\'erie
enti\`ere $\zeta(3)U(t)-V(t)$. 
En particulier, $U(t)$ et $V(t)$ sont des $G$-fonctions de Siegel. La 
s\'erie $U(t)$ est une solution de l'\'equation diff\'erentielle 
$\Lmod y=0$~; la conjecture de Bombieri-Dwork pr\'edit 
(\cite{DworkAmice}, \cite{DworkINDAM}~; voir aussi \cite{Andre}
et \cite{GerottoSullivan}) que $\Lmod$ provient de la
g\'eom\'etrie. 

Or, pour $t \in \Pun(\C) \moins \{0, 1, (\sqrt{2} \pm 1)^4, \infty\}$, 
Beukers et Peters construisent \cite{BeukersPeters}
une surface K3
$X_t$ birationnellement \'equivalente  \`a la surface projective
$S_t$ d'\'equation affine
$1-(1-xy)z-txyz(1-x)(1-y)(1-z)=0$.  Ils montrent que si $\om_t$ est l'unique
2-forme holomorphe sur $X_t$ (\`a proportionnalit\'e pr\`es), et si $\tau_t$ 
est un certain 2-cycle (constant pour la connexion de Gauss-Manin),
alors $U(t)$ est
l'int\'egrale de $\om_t$ sur $\tau_t$. En particulier 
$\Lmod y=0$ est l'\'equation de
Picard-Fuchs de cette famille de surfaces~: elle provient bien de la
g\'eom\'etrie.

\subsection{Congruences}

De nombreux auteurs ont \'etudi\'e   des propri\'et\'es de
congruence sur les nombres d'Ap\'ery $u_n$. Par exemple, Chowla, Cowles et
Cowles \cite{ChowlaCC} ont conjectur\'e $u_p \congru 5 \mod p^3$ pour tout $p
\geq 5$ premier. Cette conjecture a \'et\'e d\'emontr\'ee par plusieurs auteurs
(voir par exemple 
\cite{Gessel}, \cite{Sury}, \ldots).
De nombreuses autres congruences ont \'et\'e prouv\'ees, pour les nombres
d'Ap\'ery et certaines de leurs g\'en\'eralisations. 

\smallskip

Notons 
$\sum_{n \geq 1} \gam_n q^n = q \prod_{n \geq 1} (1-q^{2n})^4 (1-q^{4n})^4$
l'unique forme parabolique normalis\'ee
de poids 4 pour $\Gamma_0(8)$.
Pour $r \geq 1$, $m \geq 1$ impair et $p$ premier impair, on a 
la congruence suivante (qui ressemble \`a celles d'Atkin - Swinnerton-Dyer, voir
\cite{Hazewinkel} \S VI.33)~:
\begin{equation} 
u_{\frac12(mp^r-1)} - \gam_p u_{\frac12(mp^{r-1}-1)} + p^3
u_{\frac12(mp^{r-2}-1)} \congru 0 \mod p^r
\end{equation}
avec la convention $u_t = 0$ si $t \notin \Z$.
Beukers la d\'emontre \cite{Beukersanother} en utilisant la construction
modulaire du paragraphe~\ref{subsecmod}. 
On en d\'eduit 
$u_{\frac{p-1}{2}} \congru \gam_p \mod p$, congruence dont 
Beukers a conjectur\'e \cite{Beukersanother} qu'elle est
vraie modulo $p^2$. Ceci a \'et\'e prouv\'e par Ishikawa \cite{IshikawaKobe}
si $p$  ne divise pas $u_{\frac{p-1}{2}}$, puis par
Ahlgren et Ono \cite{AhlgrenOnocongr} dans le cas g\'en\'eral. Ahlgren et Ono
utilisent des s\'eries hyperg\'eom\'etriques sur $\Fp$ et la
modularit\'e de la vari\'et\'e d'\'equation 
$x + \frac{1}{x} + y + \frac{1}{y} + z + \frac{1}{z} + w + \frac{1}{w}=0$
(dont la famille de surfaces K3 consid\'er\'ee par Beukers-Peters est un
quotient~: voir  \cite{PetersStienstra}, Th\'eor\`eme 4).

Pour $r,m \geq 1$ et $p \geq 5$ premier, Beukers a d\'emontr\'e 
\cite{Beukerscongr}, de mani\`ere \'el\'ementaire, 
qu'on a $u_{mp^r-1} \congru  u_{mp^{r-1}-1}
\mod p^{3r}$.  La m\^eme congruence, mais seulement modulo $p^r$,
s'interpr\`ete en disant que $\int_0 ^T U(t) \dd t$ est (vue comme s\'erie
formelle en $T$) le logarithme d'une loi de groupe formel sur $\Z$ qui est
isomorphe \`a $\Gm$ sur $\Z$ (\cite{Beukerscongr}~; voir aussi 
l'appendice de \cite{StienstraBeukers} ou \cite{Hazewinkel}, \S VI.33). 


\mathversion{bold}
\section{Irrationalit\'e d'une infinit\'e de $\zeta(2k +1)$} 
						\label{secrivoal}
\mathversion{normal}



\subsection{\'Enonc\'e des r\'esultats} \label{subsecenoncesinfi}

Dans cette partie, on d\'emontre les r\'esultats suivants, dont le 
premier implique le
th\'eor\`eme~\ref{thinfinite}~:

\begin{theo}[\cite{RivoalCRAS}, \cite{BR}] \label{theologa}
 Pour $\aeno \geq 3$ impair, notons $\del_\aeno$ la dimension du $\Q$-espace 
vectoriel engendr\'e par  1, $\zeta(3)$, $\zeta(5)$,
\ldots, $\zeta(\aeno)$. Pour tout $\eps > 0$ il existe un entier $\aeno_0$ tel 
que pour tout $\aeno \geq \aeno_0$ impair on ait~:
$$\del_\aeno \geq \frac{1-\eps}{1+\log(2)} \log(\aeno).$$
\end{theo}

\begin{rema} Si dans le th\'eor\`eme~\ref{theologa} on remplace
$\frac{1-\eps}{1+\log(2)}$ par $\frac{1}{3}$ alors \cite{BR}
on peut prendre $\aeno_0 = 3$.
\end{rema}

\begin{theo}[\cite{BR}]
\label{centsoixante} Il existe un entier impair $\aeno$, avec $\aeno \leq
169$, tel que 1, $\zeta(3)$ et $\zeta(\aeno)$ soient lin\'eairement ind\'ependants
sur $\Q$.
\end{theo}

Ce th\'eor\`eme a \'et\'e am\'elior\'e par Zudilin \cite{Zudilincentqc}, qui
remplace 169 par 145, gr\^ace \`a un raffinement du lemme~\ref{lemmeNikishin}
ci-dessous.

\smallskip

Les deux ingr\'edients essentiels de la d\'emonstration du 
th\'eor\`eme~\ref{theologa} sont l'absence de $\zeta(2)$, $\zeta(4)$, \ldots,
$\zeta(\ell-1)$ d'une part, et la minoration en $\log(\ell)$ de la dimension
d'autre part. Seule cette deuxi\`eme id\'ee est utile pour d\'emontrer le
th\'eor\`eme suivant.

\begin{theo}[\cite{theseTanguy}] 
\label{theopolylogs} Soient $z \in \Q$, $\vert z \vert > 1$, 
et $\eps > 0$. Il existe  un entier $\aeno_0$ 
(qui d\'epend de $z$ et $\eps$) tel
que, pour tout $\aeno \geq \aeno_0$, la dimension du $\Q$-espace 
vectoriel engendr\'e par  $1, \Li_1(1/z), \Li_2(1/z), \ldots, \Li_\aeno(1/z)$
soit minor\'ee par $\frac{1-\eps}{1+\log(2)} \log(\aeno)$.
\end{theo}
En cons\'equence, pour tout nombre rationnel $z$ de 
valeur absolue sup\'erieure \`a 1 il existe
une infinit\'e
d'entiers $j$ tels que $\Li_j(1/z)$ soit irrationnel. Par ailleurs, 
quand $z$ est un entier
n\'egatif tel que $\vert z \vert > (4\aeno)^{\aeno(\aeno-1)}$, 
Nikishin a d\'emontr\'e
\cite{Nikishin} que les nombres $1, \Li_1(1/z), \Li_2(1/z), \ldots, \Li_\aeno(1/z)$
sont lin\'eairement ind\'ependants sur $\Q$~; sa m\'ethode 
a inspir\'e en partie la construction
expos\'ee au paragraphe suivant. Hata a raffin\'e 
(\cite{Hatapolylogs}, \cite{Hatadilog}) le
r\'esultat de Nikishin~: par exemple 
$1$, $\Li_1(1/z)$ et $\Li_2(1/z)$ sont lin\'eairement 
ind\'ependants sur $\Q$ pour $z \leq -5$ ou $z \geq 7$.

\subsection{Structure de la preuve} \label{subsecstruct}

Soient $a$ et $r$ deux entiers, avec $a \geq 3$ et
$1 \leq r < \frac{a}{2}$.
Soit $n \geq 1$. D\'efinissons $\RRR_n$ et $\SSS_n$ (qui d\'ependent 
aussi de $a$ et $r$) par~:
\begin{multline*}
\RRR_n(k) = 2n!^{a-2r} (k+\frac{n}{2}) \frac{(k-rn)_{rn} 
(k+n+1)_{rn}}{(k)_{n+1}^a} \\
= 2n!^{a-2r} (k+\frac{n}{2}) \frac{(k-1)(k-2)\ldots(k-rn) (k+n+1)
(k+n+2)\ldots (k+(r+1)n)}{k^a (k+1)^a \ldots (k+n)^a} 
\end{multline*}
et
\begin{equation} \label{eqdefs}
\SSS_n(z) = \sum_{k \geq 1} \RRR_n(k) z^{-k}.
\end{equation}
Cette s\'erie converge absolument pour tout nombre complexe $z$ tel que 
$\vert z \vert \geq 1$, car 
$\RRR_n(k)
= \gdo (k^{-2})$ quand $k$ tend vers l'infini.

\medskip

Les propri\'et\'es de cette s\'erie \'etudi\'ees au 
paragraphe~\ref{subsecdetails} permettent de d\'emontrer 
les th\'eor\`emes~\ref{theologa} (en prenant $z=1$ et $a$ 
pair),~\ref{centsoixante}
(avec $z=1$, $a= 169$,  $r = 10$ et $n$ impair~; on utilise
le th\'eor\`eme d'Ap\'ery) 
et~\ref{theopolylogs} (avec $z \in \Q$, $z > 1$~; pour $z < -1$
il suffirait de modifier le lemme~\ref{lemanal}).
Les trois preuves sont parall\`eles~; on d\'etaille dans ce paragraphe 
la structure de celle du th\'eor\`eme~\ref{theologa}.

\medskip

On suppose $a$ pair~; on construit des formes lin\'eaires en 
1, $\zeta(3)$, $\zeta(5)$, \ldots, $\zeta(a-1)$ gr\^ace \`a la proposition
suivante~:
\begin{prop} \label{proprivoal} Supposons $a$ pair.
Notons $d_n$ le p.p.c.m des entiers de 1 \`a $n$. 
Alors il existe des nombres rationnels $\coez$, $\coe_3$, $\coe_5$, 
\ldots, $\coe_{a-1}$ tels que~:
\begin{enumerate}
\item On a $\SSS_n(1) = \coez + \coe_3 \, \zeta(3) + 
\coe_5 \, \zeta(5) + \coe_7 \, \zeta(7) + \ldots + \coe_{a-1} \, \zeta(a-1)$.
\item \label{assertionnest} Pour tout $j \in \{0, 3, 5, \ldots, a-1\}$ on a \,
$\limsup _{n \rightarrow + \infty} \vert \coe_j \vert ^{1/n} \leq
2^{a-2r} (2r+1)^{2r+1}$.
\item \label{ameliorable}
Pour tout $j \in \{0,3,5,\ldots, a-1\}$, 
le nombre rationnel $d_n^{a} \coe_j$ est un entier.
\item Il existe un r\'eel $\psi_{r,a} > 0$ tel que \,
$\lim_{n \rightarrow + \infty} \vert \SSS_n(1)  \vert ^{1/n} = \psi_{r,a}
\leq \frac{2^{r+1}}{r^{a-2r}}$.
\end{enumerate}
\end{prop}


En fait on conjecture que l'am\'elioration suivante est possible~:

\begin{conj}[\cite{theseTanguy}] \label{conjball}
Dans l'assertion \eqref{ameliorable} de la proposition~\ref{proprivoal}, on
peut remplacer $d_n ^a$ par $d_n^{a-1}$.
\end{conj}

\begin{rema} En prenant $a=4$ (et $r=1$), on obtient les formes lin\'eaires 
en 1 et $\zeta(3)$ du paragraphe~\ref{subsecball}, donc la conjecture~\ref{conjball}
est vraie quand $a = 4$. Elle est d\'emontr\'ee aussi quand
$a=6$ et $r=1$ (voir la fin du paragraphe \ref{subsecremarques}).
On ne conna\^{\i}t pas de cons\'equence directe de cette conjecture, mais
une version forte de celle-ci pourrait \'eventuellement 
permettre de d\'emontrer que parmi $\zeta(5)$, $\zeta(7)$ et $\zeta(9)$, 
l'un au moins est irrationnel (voir la remarque  \ref{remaneuf}).
En tout cas, il serait int\'eressant d'obtenir une 
preuve  de la conjecture~\ref{conjball} gr\^ace \`a une interpr\'etation
(par exemple g\'eom\'etrique, comme au paragraphe~\ref{subsecmod})
de $\coez$, \ldots, $\coe_{a-1}$.
\end{rema}

\bigskip

La proposition~\ref{proprivoal} fournit des formes lin\'eaires en 1, $\zeta(3)$,
\ldots, $\zeta(a-1)$ (si $a$ est pair). Si cette suite de formes lin\'eaires
tend vers 0, sans \^etre nulle \`a partir d'un certain rang, 
alors l'un au moins des nombres $\zeta(3)$,
\ldots, $\zeta(a-1)$ est irrationnel. Cette remarque sera utilis\'ee 
pour d\'emontrer le th\'eor\`eme \ref{theoonze}.
Ici on veut obtenir les th\'eor\`emes~\ref{theologa} \`a~\ref{theopolylogs},
donc on a besoin d'un crit\`ere d'ind\'ependance lin\'eaire, 
qui donne une minoration plus fine
de la dimension du $\Q$-espace vectoriel engendr\'e par 1, $\zeta(3)$,
\ldots, $\zeta(a-1)$. On va utiliser \`a cet effet
le th\'eor\`eme~\ref{thnest} ci-dessous.

La meilleure minoration qu'on puisse esp\'erer 
est
donn\'ee par le principe des tiroirs, de la mani\`ere suivante. 
Soient $\alpha$  et $\beta$ des r\'eels, 
avec $0<\alpha<1$ et $\beta >1$. Soient $\theta_1,\ldots,\theta_s$
des r\'eels qui engendrent un $\Q$-espace vectoriel de dimension au moins
$1 - \frac{\log(\alpha)}{\log(\beta)}$. Alors il existe une suite $(\ell_n)$ 
de formes lin\'eaires en $\theta_1$, \ldots, $\theta_s$
dont les coefficients entiers $p_{j,n}$ v\'erifient
$\limsup_{n \rightarrow + \infty} \vert p_{j,n} \vert^{1/n} 
\leq \beta$ pour tout $j$ et telle que
$\limsup_{n \rightarrow + \infty} \vert \ell_n(\theta_1,\ldots,
\theta_s) \vert ^{1/n} \leq \alpha$. Essentiellement, plus 
la dimension du
$\Q$-espace vectoriel engendr\'e est grande, plus les formes lin\'eaires qu'on
peut construire sont petites. On cherche une
r\'eciproque \`a cette assertion.
Une contrainte suppl\'ementaire est n\'ecessaire~: si 
$\frac{\theta_2}{\theta_1}$ est un nombre de Liouville, on peut construire 
des formes lin\'eaires extr\^emement petites m\^eme si 
la dimension du 
 $\Q$-espace vectoriel engendr\'e est seulement 2. 
Ce contre-exemple ne tient plus si on demande que les formes lin\'eaires
en $\theta_1,\ldots,\theta_s$ ne soient pas trop petites. On a alors la
r\'eciproque suivante (pour une preuve, voir
\cite{Nesterenkocritere} ou \cite{Colmez}, \S II.1)~:


\begin{theo}[\cite{Nesterenkocritere}] \label{thnest} Soient 
$\theta_1,\ldots,\theta_s$ des r\'eels. Pour tout $n \geq 1$, soit
$\ell _n = p_{1,n} X_1 + \ldots + p_{s,n}
X_s$ une forme lin\'eaire \`a coefficients entiers. Soient  
$\alpha$  et $\beta$ des r\'eels, avec $0<\alpha<1$ et $\beta >1$.

Supposons qu'on ait $\limsup_{n \to + \infty} 
\vert p_{j,n} \vert^{1/n}  \leq  \beta$ pour tout $j$ compris entre 1 et $s$,
et $$\lim_{n \to + \infty} \vert \ell_n(\theta_1,\ldots,
\theta_s) \vert^{1/n} = \alpha.$$
Alors le $\Q$-espace vectoriel engendr\'e par 
$\theta_1,\ldots,\theta_s$ est de dimension au moins
$1 - \frac{\log(\alpha)}{\log(\beta)}$.
\end{theo}

Pour d\'eduire le th\'eor\`eme~\ref{theologa} de la proposition~\ref{proprivoal}
et de ce 
crit\`ere d'ind\'ependance lin\'eaire,
il suffit de consid\'erer $d_n ^a \SSS_n(1)$,
qui est une forme lin\'eaire \`a coefficients entiers
en 1, $\zeta(3)$, $\zeta(5)$,
\ldots, $\zeta(a-1)$. On choisit 
$a$ suffisamment grand, et $r$ \'egal \`a la partie enti\`ere de 
$\frac{a}{(\log(a))^2}$. Alors 
$r^r $ est n\'egligeable devant $c^a$ (pour toute constante
$c$), et on peut prendre $\beta$ essentiellement \'egal \`a 
$(2e)^a = e^{a(1+\log(2))}$ et $\alpha$ essentiellement major\'e 
par $r^{-a}$, qui est de l'ordre
de $e^{-a \log(a)}$. Cela d\'emontre le th\'eor\`eme~\ref{theologa}.


\subsection{Quelques d\'etails sur la preuve} \label{subsecdetails}

Soit $z$ un nombre complexe de module sup\'erieur ou \'egal \`a 1. 
La s\'erie $\SSS_n(z)$ peut s'\'ecrire comme une s\'erie
hyperg\'eom\'etrique tr\`es bien \'equilibr\'ee, de la mani\`ere suivante~:
\begin{multline*}
\SSS_n(z) = z^{-rn-1} n!^{a-2r}
\frac{(rn)!((r+1)n+2)_{rn+1}}{(rn+1)_{n+1} ^a} \times \\
\fpqusz{a+3}{a+2}{(2r+1)n+2}{(r+\frac{1}{2})n+2}{rn+1}{rn+1}{(r+
\frac{1}{2})n+1}{(r+1)n+2}{(r+1)n+2}.
\end{multline*}
Cette identit\'e provient  de simplifications dans les symboles
de Pochhammer.

\subsubsection{Repr\'esentation int\'egrale et estimation analytique}

On a la repr\'esentation int\'egrale suivante, pour $\vert z \vert \geq 1$~:
\begin{small}
$$\SSS_n(z) = 
\frac{((2r+1)n+2)!}{n!^{2r+1}} z^{(r+1)n+1}
\int_{[0,1]^{a+1}} \left( \frac{\prod_{j=1} ^{a+1} t_j ^r (1-t_j)}{(z-
t_1 t_2\ldots t_{a+1})^{2r+1}}\right) ^n \frac{z+t_1\ldots t_{a+1}}{(z-
t_1\ldots t_{a+1})^3} \dd t_1 \ldots \dd t_{a+1}.$$ 
\end{small}
Cette formule (voir par exemple \cite{Catalan}, Lemme 1) se d\'eduit 
de l'\'ecriture de $\SSS_n(z)$ comme s\'erie hyperg\'eom\'etrique~: pour
$\vert z \vert > 1$ on applique les 
relations (4.1.2) et (1.5.21) de \cite{Slater}, puis on prolonge \`a 
$\vert z \vert = 1$ par continuit\'e (voir la preuve du lemme~2 de \cite{BR}).
On peut aussi obtenir une preuve directe en d\'eveloppant en s\'erie le
d\'enominateur de l'int\'egrande (\cite{Colmez}, \cite{Habsieger}).

En calculant le maximum
sur $[0,1]^{a+1}$ de la fonction dont on int\`egre la puissance $n$-i\`eme,
on d\'eduit de cette repr\'esentation int\'egrale l'estimation analytique
suivante~:

\begin{lemm} \label{lemanal}
On suppose $z \in \R$, $z \geq 1$.
Le polyn\^ome
$$Q_{r,a,z}(s) = r s^{a+2} - (r+1)s^{a+1}+(r+1)zs - rz$$
admet une racine unique $s_0 \in [0,1]$, et elle v\'erifie $s_0 >
\frac{r}{r+1}$. De plus, si
$$\phi_{r,a,z} = z^{-r}  ((r+1)s_0 - r)^r (r+1-rs_0)^{r+1}(1-s_0)^{a-2r},$$
alors
$$\lim_{n \to \infty} \vert \SSS_n(z) \vert ^{1/n} = \phi_{r,a,z}
\leq \frac{2^{r+1}}{z^r r^{a-2r}}.$$
\end{lemm} 


Pour d\'emontrer ce lemme, il 
suffit d'adapter les preuves du lemme 2.2 de \cite{theseTanguy} et du
lemme 3 de \cite{BR}. On pourrait aussi
donner une d\'emonstration \'el\'ementaire 
de ce comportement asymptotique, sans utiliser la repr\'esentation int\'egrale
(comme la deuxi\`eme preuve du lemme 3 de \cite{BR}). Enfin, une troisi\`eme
possibilit\'e serait d'\'ecrire $\SSS_n(z)$ comme int\'egrale complexe 
et d'appliquer la m\'ethode du col~; mais cette m\'ethode est tr\`es difficile
\`a mettre en \oe uvre quand $r$, $a$ et $z$ sont des param\`etres.

\begin{rema} Pour d\'emontrer les th\'eor\`emes~\ref{theologa} et~\ref{theopolylogs}, il suffit de conna\^{\i}tre
l'existence de la limite de $\vert \SSS_n(z) \vert ^{1/n}$,
et  sa majoration par
$\frac{2^{r+1}}{z^r r^{a-2r}}$. La valeur exacte de $\phi_{r,a,z}$ n'est utile
que pour obtenir des estimations num\'eriques pr\'ecises (par
exemple pour le th\'eor\`eme~\ref{centsoixante}).
\end{rema}

\subsubsection{D\'ecomposition en polylogarithmes} \label{subsubsecdcp}

Pour d\'emontrer que $\SSS_n(z)$ est une combinaison lin\'eaire (\`a coefficients
rationnels) de $1$, $\Li_1(1/z)$, \ldots, $\Li_a(1/z)$ quand 
$\vert z \vert > 1$, il suffit de d\'ecomposer 
la fraction rationnelle $\RRR_n$ en \'el\'ements simples,
sous la forme suivante~:
\begin{equation} \label{dcpelsples}
\RRR_n (k) = \sum_{i=0} ^n \sum_{j=1} ^a \frac{\cij}{(k+i)^j}
\end{equation}
o\`u les coefficients $\cij$ sont des rationnels, donn\'es par
\begin{equation} \label{eqdefcij}
\cij = \frac{1}{(a-j)!} \left( \frac{\dd}{\dd X} \right) ^{a-j} 
(\RRR_n(X) (X+i)^a) _{\vert X = -i}.
\end{equation}
On a pour $\vert z \vert > 1$~:
\begin{eqnarray*}
\SSS_n(z) &=&   \sum_{i=0} ^n \sum_{j=1} ^a \cij
 \sum_{k \geq 1} \frac{z^{-k}}{(k+i)^j} \\
 &=&  \sum_{i=0} ^n \sum_{j=1} ^a \cij z^i  \Li_j(1/z) -
  \sum_{i=0} ^n \sum_{j=1} ^a \cij \sum_{q=1} ^i \frac{z^{i-q}}{q^j},
\end{eqnarray*}
d'o\`u
\begin{equation} \label{eqdcppolylogs}
\SSS_n(z) = P_0(z) + \sum_{j=1} ^a P_j(z) \Li_j(1/z)
\end{equation} 
en posant
\begin{equation} \label{eqdefpzero}
P_0(z) = - \sum_{\ell = 0} ^{n-1} \left( \sum_{i=\ell+1} ^n \sum_{j=1} ^a 
\frac{\cij}{(i-\ell)^j} \right) z^\ell
\end{equation} 
et
\begin{equation} \label{eqdefpj}
P_j(z) = \sum_{i=0} ^n \cij z^i \mbox{ pour } j \in \{1, \ldots,a\}.
\end{equation} 
Bien s\^ur, les $P_j$ et les $\cij$ d\'ependent aussi de $n$, $a$ et $r$.

\subsubsection{Propri\'et\'e de sym\'etrie} \label{subsubsecsym}

La fonction $\RRR_n$ v\'erifie la propri\'et\'e 
de sym\'etrie suivante~:
$$\RRR_n(-k-n) = (-1)^{a(n+1)+1} \RRR_n(k).$$
Cette sym\'etrie est rendue possible par la pr\'esence des deux facteurs de
Pochhammer au num\'erateur de $\RRR_n(k)$~: quand $k$ est
chang\'e en $-k-n$, ils sont permut\'es 
(on applique la formule $(-\al)_p = (-1)^p (\al - p +1)_p$).

L'unicit\'e du d\'eveloppement en \'el\'ements simples montre que
$c_{i,j} = (-1)^{j+a(n+1)+1}c_{n-i,j}$
pour tous $i \in \zeron$ et $j \in \una$,
ce qui donne pour tout $j  \in \una$~:
\begin{equation} \label{recipropol}
P_j(z) = (-1)^{j+a(n+1)+1} z^n P_j(1/z).
\end{equation}
En particulier, si $j+a(n+1)$ est pair alors $P_j(1) = 0$. De plus on a 
$P_1(1) = 0$, car
$P_1(1) = \sum_{i=0} ^n c_{i,1}$ est l'oppos\'e du r\'esidu \`a l'infini de
$\RRR_n$ (on peut aussi faire tendre $z$ vers 1 dans \eqref{eqdcppolylogs}
 et constater que le seul terme qui puisse tendre vers l'infini est $P_1(z)
 \Li_1(1/z)$). Quand $a$ est pair,
on obtient donc~:
$$\SSS_n(1) = P_0(1) + P_3(1) \zeta(3) +  P_5(1) \zeta(5) +  \ldots
+ P_{a-1}(1) \zeta(a-1).$$ 
Quand $a$ est impair et $n$ pair, on obtient de m\^eme une forme 
lin\'eaire en 
1, $\zeta(2)$, $\zeta(4)$, \ldots, $\zeta(a-1)$ dont on peut se servir pour 
montrer
qu'une infinit\'e de puissances de $\pi$ sont lin\'eairement ind\'ependantes sur
$\Q$, i.e. que $\pi$ est transcendant. On peut aussi 
en d\'eduire une mesure de 
transcendance de $\pi$, \`a la mani\`ere de Reyssat
\cite{Reyssatmesure}. 

Enfin, quand $a$ et $n$ sont impairs, on obtient une forme lin\'eaire en
1, $\zeta(3)$, $\zeta(5)$, \ldots, $\zeta(a)$~; c'est ce qu'on utilise pour
d\'emontrer le th\'eor\`eme~\ref{centsoixante}.

\subsubsection{Majoration des coefficients de la forme lin\'eaire}
						\label{subsecmajcoeffs}
\begin{lemm} Pour tout $j \in \zeroa$ on a~:
$$\limsup _{n \rightarrow + \infty} \vert P_j(z) \vert ^{1/n} \leq
2^{a-2r} (2r+1)^{2r+1}\vert z \vert .$$
\end{lemm}

\Dem On peut suivre la d\'emonstration 
du lemme 4 de \cite{BR} en \'ecrivant  la
formule de Cauchy sur le cercle $\cercle$ de centre $-i$ et de rayon
$1/2$~:
$$c_{i,j} = \frac{1}{2 i \pi} \int_{\cercle} \RRR_n(t) (t+i)^{j-1} \dd t.$$
On majore ensuite le module de 
l'int\'egrande, et le lemme en d\'ecoule.
Une autre preuve, qui conduit \`a une majoration l\'eg\`erement moins pr\'ecise,
est donn\'ee dans 
\cite{Colmez} et \cite{Habsieger}.

\subsubsection{Estimation arithm\'etique}

Les polyn\^omes $P_0,\ldots,P_a$ sont \`a coefficients rationnels~; on a besoin
d'un d\'enominateur commun pour leurs coefficients.

\begin{lemm} \label{lemmeNikishin}
Pour
tout $j \in \{0,\ldots,a\}$, le polyn\^ome $d_n ^{a-j} P_j(z)$ est \`a 
coefficients entiers.
\end{lemm}

\begin{rema} On peut (\cite{Zudilincentqc}, \S 4) raffiner ce lemme, 
ce qui permet
de remplacer 169 par 145 dans l'\'enonc\'e du th\'eor\`eme~\ref{centsoixante}.
Cependant, des exemples montrent qu'on ne peut pas esp\'erer remplacer 
$d_n ^{a-j}$ par $d_n ^{a-1-j}$. La conjecture~\ref{conjball} signifie que pour $z=1$ on a des compensations particuli\`eres
qui font chuter le d\'enominateur. 
\end{rema}

\Dem 
Posons $F_s (X) = \frac{(X-sn)_n}{(X)_{n+1}}$ et 
$G_s (X) = \frac{(X+sn+1)_n}{(X)_{n+1}}$ pour tout $s \in \{1,\ldots, r\}$,
ainsi que
$H(X)   = \frac{n!}{(X)_{n+1}}$ et 
$I(X) = 2X+n$.
Alors on a $F_s(X) = \sum_{p=0}^{n} \frac{f_{p,s}}{X+p}$ avec
$f_{p,s}=(-1)^{n-p} \combi{n}{p} \combi{p+sn}{n} \in \Z$, et de m\^eme 
(avec des notations \'evidentes) $g_{p,s} \in \Z$ et $h_p \in \Z$ pour tous
$p$, $s$. On obtient alors le d\'eveloppement en \'el\'ements simples de 
$\RRR_n(X) = \left( \prod_{s=1} ^r F_s (X) \right) \cdot
\left( \prod_{s=1} ^r G_s (X) \right) \cdot H(X)^{a-2r} \cdot I(X)$
en faisant le produit des d\'eveloppements des facteurs. 
On utilise les formules $\frac{2X+n}{X+p} = 2 + \frac{n-2p}{X+p}$ et  
$\frac{1}{(X+p)(X+p')} = \frac{1}{(p'- p)(X+p)} + \frac{1}{(p- p')(X+p')}$
pour $p \neq p'$~; les d\'enominateurs n'apparaissent que par application de la
seconde.
Ce calcul montre que $d_n ^{a-j} \cij$ est entier pour tous $i$, $j$, ce qui
ach\`eve la preuve (suivant \cite{Colmez} et 
\cite{Habsieger}) du lemme.


\subsection{Quelques remarques} \label{subsecremarques}

Soit $Q_n$ un polyn\^ome \`a coefficients rationnels, de degr\'e inf\'erieur ou
\'egal \`a $a(n+1)-1$. On peut toujours
consid\'erer $\RRR_n(k) = \frac{Q_n(k)}{(k)_{n+1}^a}$ et
$\SSS_n(z) = \sum_{k \geq 1} \RRR_n(k) z^{-k}$, qui converge 
quand $\vert z \vert > 1$. Une difficult\'e majeure consiste \`a bien choisir le
polyn\^ome $Q_n$.

Quel que soit ce choix, on peut d\'ecomposer $\RRR_n$ en \'el\'ements simples,
d\'efinir $P_0$, \ldots, $P_a$ et obtenir une d\'ecomposition de $\SSS_n(z)$ en
polylogarithmes~: toutes les formules du paragraphe~\ref{subsubsecdcp} restent
valables. Pour obtenir une forme lin\'eaire en valeurs de $\zeta$, il 
faut\footnote{Voir cependant la remarque \ref{remamoinsun}.}
faire tendre $z$ vers 1. 
Tous les termes de la d\'ecomposition en
polylogarithmes ont une limite finie, sauf peut-\^etre $P_1(z) \Li_1(1/z)$.
C'est pourquoi on suppose $P_1(1)=0$, ce qui signifie que $\RRR_n$ n'a pas de
r\'esidu \`a l'infini, i.e. $\deg(Q_n) \leq a(n+1)-2$~; alors 
la s\'erie qui d\'efinit $\SSS_n(z)$ converge absolument d\`es que $\vert z
\vert \geq 1$.

En outre on souhaite\footnote{Sauf pour d\'emontrer le th\'eor\`eme
\ref{theopolylogs}~; pour ce dernier, le polyn\^ome $Q_n(k) = (k-rn)_{rn}$
convient aussi. C'est celui qui est utilis\'e dans le Chapitre 2 de 
\cite{theseTanguy}.} 
obtenir une forme lin\'eaire en les
$\zeta(2k+1)$ seulement, c'est-\`a-dire avoir $P_{j}(1)=0$ pour tout $j \geq 2$
pair. Pour assurer cela il est suffisant d'avoir une propri\'et\'e de
sym\'etrie du polyn\^ome $Q_n$, en l'occurrence 
$Q_n(-k-n) = (-1)^{a(n+1)+1} Q_n(k)$. 
C'est cette remarque qui constitue le c\oe ur des progr\`es r\'ecents
(\cite{RivoalCRAS}, \cite{BR}). On ne sait pas du tout la g\'en\'eraliser, par
exemple pour  construire 
des formes lin\'eaires en $\zeta(s)$ dans lesquelles les $s$
appartenant \`a une certaine progression arithm\'etique n'apparaissent pas.

La forme lin\'eaire $\SSS_n(1)$ ne sera int\'eressante que
si elle tend suffisamment vite vers 0 quand $n$ tend vers l'infini. 
Intuitivement, ce sera le cas si les premiers termes de la s\'erie qui d\'efinit
$\SSS_n(1)$ sont nuls. C'est pourquoi on cherche un polyn\^ome $Q_n(k)$ qui
s'annule aux premiers entiers, en l'occurrence entre  1 et $rn$~; ceci
signifie que $Q_n(k)$ est multiple de $(k-rn)_{rn}$. Il s'agit en fait d'un
probl\`eme de type Pad\'e~: on demande aux polyn\^omes $P_0$, \ldots, $P_a$
d'\^etre tels que
$$\SSS_n(z) =
P_0(z) + \sum_{j=1} ^a P_j(z) \Li_j(1/z) = \gdo (z^{-rn-1}) \mbox{ quand }
z \to \infty.$$

Parmi tous les polyn\^omes sym\'etriques $Q_n(k)$ multiples 
de $(k-rn)_{rn}$ (donc n\'ecessairement aussi multiples de $(k+n+1)_{rn}$),
on a int\'er\^et \`a en prendre un de degr\'e minimal, pour que $\SSS_n(1)$ soit
aussi petit que possible. Si $a(n+1)$ est impair, le polyn\^ome
$(k-rn)_{rn}(k+n+1)_{rn}$ a la bonne parit\'e, et on peut consid\'erer 
$Q_n(k) = n!^{a-2r} (k-rn)_{rn}(k+n+1)_{rn}$~: on obtient la s\'erie
hyperg\'eom\'etrique bien \'equilibr\'ee de \cite{RivoalCRAS}
et \cite{BR}. 
Si $a(n+1)$ est pair, pour obtenir le bon signe dans la
propri\'et\'e de sym\'etrie de $Q_n$ on est amen\'e \`a introduire
un facteur $k + \frac{n}{2}$, ce qui donne la s\'erie tr\`es bien
\'equilibr\'ee du paragraphe~\ref{subsecstruct}. Dans les deux cas, 
$\SSS_n(z)$ est la solution unique d'un probl\`eme de
Pad\'e (voir \cite{Huttner03} et \cite{FischlerRivoal}).

Plus $a$ est grand (en prenant, pour chaque $a$, la valeur optimale de $r$), 
plus la forme lin\'eaire \`a
coefficients entiers $d_n ^a \SSS_n(1)$ 
est petite (et la pr\'esence, ou l'absence, du facteur $k
+ \frac{n}{2}$ a une influence n\'egligeable sur ce comportement). Donc 
si on cherche des formes lin\'eaires en 1, $\zeta(3)$, $\zeta(5)$, \ldots,
$\zeta(2\ell+1)$, celles obtenues avec
la s\'erie tr\`es bien \'equilibr\'ee
pour $a = 2\ell+2$ seront meilleures que celles obtenues avec
la s\'erie bien
\'equilibr\'ee pour $a = 2\ell+1$ et $n$ pair. Ceci n'a aucune 
influence quand $\ell$
tend vers l'infini, mais peut s'av\'erer crucial si $\ell$ est fix\'e 
(comme dans le th\'eor\`eme~\ref{theoonze}). 
En outre, si la conjecture~\ref{conjball} 
(qui n'a aucun \'equivalent pour des
s\'eries seulement bien \'equilibr\'ees)
est vraie alors il suffit de multiplier $\SSS_n(1)$ par
$d_n ^{a-1}$, ce qui donne une forme lin\'eaire encore plus petite. 
Pour $a=4$, on retrouve ainsi les formes lin\'eaires 
d'Ap\'ery en 1 et $\zeta(3)$ (ce qui n'est pas le cas avec la s\'erie bien
\'equilibr\'ee quand $a=3$).

\begin{rema} \label{remamoinsun} Pour d\'emontrer le th\'eor\`eme
\ref{theologa} on pourrait \'evaluer les formes lin\'eaires en
polylogarithmes en $z = -1$  plut\^ot qu'en $z=1$. Ceci induit peu de
changements. Le plus notable est que $\log(2) = -\Li_1(-1)$
remplace le divergent $\Li_1(1)$~; pour $\ell \geq 2$ on a 
$\Li_\ell (-1) = -(1-2^{1-\ell})\zeta(\ell)$. 
Pour $a=3$ et $z=-1$ les formes lin\'eaires 
construites au paragraphe \ref{subsecdetails} sont \cite{Krattenthaler}
celles utilis\'ees par Ap\'ery (\cite{Apery}, \cite{VDP})
pour prouver que $\zeta(2)$
est irrationnel. En particulier $d_n^2$ suffit comme d\'enominateur des 
coefficients de cette forme lin\'eaire. Plus g\'en\'eralement, la 
conjecture~\ref{conjball} devrait \^etre valable aussi quand $a$ est 
impair et $z = -1$.
\end{rema}


\bigskip

Consid\'erons l'op\'erateur diff\'erentiel hyperg\'eom\'etrique 
suivant, o\`u $\de = z\frac{\dd}{\dd z}$~:
$$\rivo = \de^{a+1} (\de - \frac{n}2 -1) (\de - (r+1)n - 1) - z
(\de-n)^{a+1} (\de-\frac{n}2+1)(\de+rn+1).$$
L'\'ecriture de $\SSS_n(z)$ 
comme s\'erie hyperg\'eom\'etrique tr\`es bien \'equilibr\'ee montre que
$\SSS_n(z)$ est une solution de l'\'equation diff\'erentielle
$\rivo y =0$. Par monodromie 
on voit, gr\^ace \`a 
\eqref{eqdcppolylogs}, que pour tout $b \in \{1, \ldots, a\}$ la fonction
$\sum_{j=b} ^a (-1)^{j-1} P_j(z) \frac{\log^{j-b}(z)}{(j-b)!}$ est aussi une
solution de $\rivo y =0$. En particulier pour $b=a$ on obtient le polyn\^ome
$P_a$, qu'on peut \'ecrire comme polyn\^ome hyperg\'eom\'etrique 
tr\`es bien \'equilibr\'e (avec un
petit abus de langage~: ici 
les param\`etres inf\'erieurs  $-\frac{n}{2}$ et $-(r+1)n$
sont n\'egatifs, mais
la s\'erie $_{a+3} F_{a+2}$ est quand m\^eme bien d\'efinie)~:
\begin{eqnarray*}
P_a(z) &=& (-1)^{rn} n (rn)! ((r+1)n)! n!^{-2r-1} \croix \\
 && _{a+3} F_{a+2} \left(
\begin{array}{cccccc|}
 -n , & -\frac{n}{2}+1 ,& rn + 1,  & -n ,& \ldots,  & -n     \\
     &  -\frac{n}{2},  & -(r+1)n ,& 1 , &  \ldots ,& 1 
\end{array}
\, \, \, z \right).
\end{eqnarray*}
L'aspect bien \'equilibr\'e de ce polyn\^ome hyperg\'eom\'etrique 
lui conf\`ere (voir  \cite{Andrews} ou \cite{AAR}, \S 3.5) 
la propri\'et\'e de r\'eciprocit\'e \eqref{recipropol}. 
En effet, si $y(z)$ est une solution de l'\'equation diff\'erentielle $\rivo
y=0$
alors $z^n y(1/z)$ est aussi une solution de cette m\^eme \'equation.
Quant aux autres polyn\^omes $P_{a-1}$, \ldots, $P_1$, ils 
s'obtiennent par la m\'ethode de 
Frobenius (voir \cite{Ince}) et  v\'erifient, eux aussi, \eqref{recipropol}.
Toutes ces consid\'erations valent aussi pour la s\'erie bien \'equilibr\'ee 
de \cite{RivoalCRAS} et \cite{BR}, et permettent \cite{Huttner03}
d'\'ecrire celle-ci comme solution unique d'un probl\`eme de Pad\'e.



\medskip

Un autre int\'er\^et des d\'efinitions utilis\'ees dans ce texte est
que $\SSS_n(1)$ poss\`ede (pour $a$ pair)
plusieurs repr\'esentations int\'egrales 
assez simples. Tout
d'abord, on a (\cite{Zudilinservice}, Th\'eor\`eme 5)
l'int\'egrale suivante, qui g\'en\'eralise $\Ibeu(1)$ et 
les int\'egrales introduites par 
Vasilenko \cite{Vasilenko} et Vasilyev (\cite{Vasilyevancien}, 
\cite{Vasilyev})~:
\begin{equation} \label{intvasi}
\SSS_n(1) = \frac{(rn)!^2}{n!^{2r}} \int_{[0,1]^{a-1}}
\frac{\prod_{j=1} ^{a-1} x_j ^{rn} (1-x_j)^n}{(\QQ_{a-1}(x_1,\ldots,
x_{a-1}))^{rn+1}} \dd x_1 \ldots \dd x_{a-1},
\end{equation}
en posant $\QQ_{a-1}(x_1,\ldots,x_{a-1}) = 1- x_1 (1- x_2 (
 \ldots (1-x_{a-1})\ldots ))$. 
Vasilyev a d\'emontr\'e \cite{Vasilyev} que si $a=6$ et $r=1$ alors
cette int\'egrale 
s'\'ecrit $\coeprz + \coepr_3 \zeta(3) + \coepr_5 \zeta(5)$ avec 
$d_n ^5 \coeprz$, 
$d_n ^5  \coepr_3$ et $d_n ^5  \coepr_5$ entiers. Ceci prouve la
conjecture \ref{conjball} dans ce cas. Il n'est pas \'evident
que $\coeprz$, $\coepr_3$ et $\coepr_5$ soient les $P_0(1)$, $P_3(1)$ et 
$P_5(1)$ du paragraphe
\ref{subsecdetails}, mais cela d\'ecoule de l'ind\'ependance lin\'eaire
conjecturale de $1$, $\zeta(3)$ et $\zeta(5)$. 

D'autre part, en appliquant \`a \eqref{intvasi}
un th\'eor\`eme de  Zlobin \cite{Zlobin} ou 
le changement de variables qui figure dans \cite{SFCRAS} (\S 2)
on obtient
l'int\'egrale suivante, qui ressemble \`a celles utilis\'ees par
Sorokin (\cite{Sorokinpi}, \cite{SorokinApery})~:
\begin{small}
$$\SSS_n(1) = \frac{(rn)!^2}{n!^{2r}} \int_{[0,1]^{a-1}}
\frac{\prod_{j=1} ^{a-1} x_j ^{rn} (1-x_j)^n \dd x_j}{(1-x_1x_2)^{n+1}
(1-x_1x_2x_3x_4)^{n+1} \ldots (1-x_1\ldots x_{a-2})^{n+1}
(1-x_1\ldots x_{a-1})^{rn+1}}.$$
\end{small}
Il serait int\'eressant d'arriver \`a d\'emontrer le 
th\'eor\`eme~\ref{theologa} en utilisant seulement des 
int\'egrales multiples comme celle-ci
(ou celle de \eqref{intvasi}). Le
probl\`eme est qu'a priori on s'attend \`a ce qu'une telle int\'egrale 
$(a-1)$-uple soit une forme
lin\'eaire, \`a coefficients rationnels, en les polyz\^etas de poids au plus
$(a-1)$ 
(voir \cite{MiW} et \cite{Zlobin}, Th\'eor\`eme~3). 
Or le th\'eor\`eme~5 de \cite{Zudilinservice} montre que ces int\'egrales sont
\'egales \`a $\SSS_n(1)$, donc seuls 1 et les valeurs de $\zeta$ aux 
entiers impairs apparaissent. 




\section{R\'esultats quantitatifs}

\mathversion{bold}
\subsection{Exposant d'irrationalit\'e de $\zeta(3)$}
\mathversion{normal}

On appelle {\em exposant d'irrationalit\'e} d'un nombre r\'eel irrationnel
$\alpha$, et on note $\mu(\alpha)$, la borne inf\'erieure de l'ensemble des
r\'eels $\nu$ pour lesquels il n'existe qu'un nombre fini de nombres rationnels 
$p/q$ tels que $\vert \alpha - \frac{p}{q} \vert < \frac{1}{q^\nu}$.
La th\'eorie des fractions continues (\cite{HW}, \S 11.1), ou
le principe des tiroirs de 
Dirichlet (\cite{HW}, \S 11.3), montre qu'un exposant 
d'irrationalit\'e est toujours sup\'erieur ou \'egal \`a 2.
Si $\alpha$ est alg\'ebrique, 
Liouville a d\'emontr\'e (\cite{Liouville} ; voir aussi
\cite{HW}, \S 11.7) que  $\mu(\alpha)$ est inf\'erieur ou \'egal au
degr\'e de $\alpha$. Ce r\'esultat a \'et\'e am\'elior\'e par Roth en 1955~: 
on a $\mu(\al) = 2$ 
pour tout nombre alg\'ebrique irrationnel $\al$ (voir \cite{EMS}, Chapitre 1, \S
7). On a aussi $\mu(\al) = 2$ pour presque tout r\'eel $\al$, au sens de la
mesure de Lebesgue (\cite{HW}, \S 11.11).
\`A l'oppos\'e, un nombre de Liouville est un 
nombre dont l'exposant d'irrationalit\'e est infini~: il est extr\^emement bien
approch\'e par des nombres rationnels (un exemple de tel nombre est 
$\sum_{k \geq 1} \frac{1}{10^{k!}}$).

\medskip

Les formes lin\'eaires d'Ap\'ery montrent que l'exposant d'irrationalit\'e 
de $\zeta(3)$ est major\'e par $13,4179$ (voir \cite{EMS}, Chapitre 2, \S 5.6)~;
en particulier $\zeta(3)$ n'est pas  un nombre de Liouville. 
Ce r\'esultat a \'et\'e am\'elior\'e notamment
par Hata \cite{Hata}
puis Rhin-Viola, qui ont d\'emontr\'e la meilleure majoration de 
$\mu(\zeta(3))$ connue \`a ce jour~:
\begin{theo}[\cite{RV3}] \label{theoRV}
L'exposant d'irrationalit\'e de $\zeta(3)$ est major\'e par $5,5139$, 
c'est-\`a-dire qu'il n'existe qu'un nombre fini de nombres
rationnels $p/q$ tels
que $$\vert \zeta(3) - \frac{p}{q} \vert < \frac{1}{q^{5,5139}}.$$
\end{theo}
Pour obtenir ce r\'esultat, Rhin et Viola consid\`erent les int\'egrales
suivantes~:
\begin{equation} \label{intrv}
\Irv = \int_0 ^1 \int_0 ^1 \int_0 ^1 
\frac{u^{hn}(1-u)^{ln}v^{kn}(1-v)^{sn}w^{jn}(1-w)^{qn}}{(1-w(1-uv))^{(q+h-r)n+1}}
\dd u \, \dd v \, \dd w,
\end{equation} 
o\`u $h, \ldots, s$ sont des param\`etres dont on fixe les valeurs de la
mani\`ere suivante :
$h=16$, $j = 17$, $k =19$, $l =15$, $q=11$, $r=9$, $s=13$.
Si on prenait tous ces param\`etres \'egaux \`a un m\^eme entier, on obtiendrait
les int\'egrales du paragraphe~\ref{subsecbeukers}, donc la suite
des formes lin\'eaires d'Ap\'ery (ou, plus pr\'ecis\'ement,
une suite extraite), conduisant \`a la m\^eme mesure d'irrationalit\'e.
L'int\'er\^et r\'eside donc dans le fait de ne pas prendre tous les param\`etres
\'egaux~; l'asymptotique obtenue pour $\Irv^{1/n}$ est un peu moins bonne, 
mais on gagne beaucoup sur les d\'enominateurs par lesquels il faut multiplier
$\Irv$ pour obtenir une forme lin\'eaire en 1 et $\zeta(3)$ \`a coefficients
entiers. Ce gain provient de l'action sur des int\'egrales de la forme
\eqref{intrv} d'un groupe isomorphe au produit
semi-direct $H \psd \scinq$, o\`u $H$ est l'hyperplan d'\'equation
$\eps_1+\ldots+\eps_5=0$ dans $(\Z / 2 \Z)^5$.
D'autres interpr\'etations de cette action de
groupe se trouvent dans 
\cite{Zudilincinqaout} et \cite{SFCaen}.


\begin{rema} Les majorations de $\mu(\zeta(3))$ mentionn\'ees ci-dessus sont
effectives~: on peut donner
une majoration explicite de la hauteur $\max(\vert p \vert, \vert q \vert)$ 
des approximations rationnelles $p/q$ ``exceptionnellement bonnes''. 
Ceci contraste avec le
th\'eor\`eme de Roth, dans lequel on sait
seulement majorer le nombre d'exceptions $p/q$, mais pas leur hauteur.
\end{rema}

\mathversion{bold}
\subsection{Irrationalit\'e d'un nombre parmi $\zeta(5)$, \ldots, $\zeta(21)$}
\mathversion{normal}

Soit $a$ un entier pair, avec $a \geq 6$.
Dans ce paragraphe, on construit (en suivant \cite{vingtetun})
des formes lin\'eaires \`a coefficients rationnels en  1,
$\zeta(5)$, $\zeta(7)$, \ldots, $\zeta(a+1)$. Si, apr\`es multiplication par un
d\'enominateur commun des coefficients, elles tendent vers z\'ero 
sans \^etre nulles \`a partir d'un certain rang, alors
l'un au moins des nombres $\zeta(5)$, $\zeta(7)$, \ldots, $\zeta(a+1)$
est irrationnel~; c'est ce qui va se produire avec $a = 20$. On pose~:
$$\Rvetun(k) = n!^{a-6}
(k+\frac{n}{2}) \frac{(k-n)_n ^3 (k+n+1)_n ^3}{(k)_{n+1} ^a}$$
et
$$\Svetun(z) = \frac12 \sum_{k=1} ^\infty \Rvetun ''(k) z^{-k}.$$ 
On d\'eveloppe $\Rvetun$ en \'el\'ements simples, ce qui d\'efinit des
coefficients $\cijvetun$ (les formules \eqref{dcpelsples} et \eqref{eqdefcij} 
restant valables). 
On d\'efinit $\Pvetun_1, \ldots, \Pvetun_a$ \`a partir des $\cijvetun$ par la relation
\eqref{eqdefpj}~; seul $\Pvetun_0$ est d\'efini par une formule l\'eg\`erement
diff\'erente~:
$$\Pvetun_0(z) = - \sum_{\ell = 0} ^{n-1} \left( \sum_{i=\ell+1} ^n \sum_{j=1} ^a 
\frac{j(j+1)\cijvetun}{2(i-\ell)^{j+2}} \right) z^\ell.$$
On obtient la
d\'ecomposition suivante  exactement comme au paragraphe~\ref{subsubsecdcp}, 
mais un d\'ecalage se produit car on d\'erive $\Rvetun$ (voir le 
paragraphe~\ref{subsecnesterenko})~:
$$\Svetun(z) = \Pvetun_0(z) + \sum_{j=1} ^a \frac{j(j+1)}{2} \, 
\, \Pvetun_j (z)
\Li_{j+2} (1/z).$$
Les arguments du paragraphe~\ref{subsubsecsym} restent valables, et montrent
(car $a$ est pair) que $\Svetun(1)$ est une forme 
lin\'eaire \`a coefficients rationnels en  1,
$\zeta(5)$, $\zeta(7)$, \ldots, $\zeta(a+1)$. De plus un d\'enominateur commun
pour ces coefficients est $2 d_n ^{a+2}$~; 
on conjecture (\cite{theseTanguy}, \S 5.1) 
que $2 d_n ^{a+1}$ convient aussi. La majoration de ces
coefficients (qui est effectu\'ee au 
paragraphe~\ref{subsecmajcoeffs}) est inutile ici~: elle servait \`a appliquer le crit\`ere
de Nesterenko, dont on n'a pas besoin puisqu'on applique seulement la remarque
\'evidente qu'une forme lin\'eaire, \`a coefficients entiers, en des
rationnels fix\'es ne peut pas \^etre arbitrairement petite sans \^etre nulle.

Le point d\'elicat de la preuve est l'estimation asymptotique de 
$\Svetun(1)$. En effet, on ne conna\^{\i}t pas d'\'ecriture de 
$\Svetun(1)$ comme int\'egrale multiple r\'eelle. 
On utilise donc la m\'ethode du col. Posons
$$\Iriv(u) = \frac{-1}{2i\pi} \int_{c- i \infty} ^{c + i \infty}
\Rvetun(s) \left(\frac{\pi}{\sin(\pi s)} \right)^3  e^{us} \dd s,$$
o\`u $c$ est un r\'eel avec $0 < c < n+1$, 
et $u$ un nombre complexe tel que $\Reelle(u) \leq 0$ et  
$\vert \Imagi(u) \vert < 3 \pi$.
Cette int\'egrale est \`a rapprocher de celle not\'ee
$\Icompl(z)$ au paragraphe~\ref{subseccplx}. On peut appliquer le th\'eor\`eme
des r\'esidus, pour faire appara\^{\i}tre les p\^oles de l'int\'egrande qui
sont situ\'es aux entiers $n+1$, $n+2$, \ldots Au voisinage d'un tel entier $k$,
on a $(\frac{\pi}{\sin(\pi s)})^3 = \frac{(-1)^k}{(s-k)^3} + \frac{(-1)^k
\pi^2}{2(s-k)} + {\small \gdo(s-k)}$. On obtient donc (voir \cite{Hessami} 
et \cite{Zudilincentqc} pour des r\'esultats analogues)~:
$$\Iriv(u) = \frac{\pi^2 + u^2}{2}
\sum_{k = n+1} ^\infty \Rvetun(k) (-e^u)^k +
u  \sum_{k = n+1} ^\infty \Rvetun ' (k) (-e^u)^k +
\frac12 \sum_{k = n+1} ^\infty \Rvetun '' (k) (-e^u)^k .$$
En choisissant $u=i \pi$, le premier terme dispara\^{\i}t, et on obtient $\Svetun(1)
= \Reelle (\Iriv(i \pi))$. 

La m\'ethode du col donne (\cite{vingtetun}, Lemme 5) deux nombres complexes 
non nuls
$c_0$ et $\alpha$, qu'on peut calculer,  
tels que $\Iriv(i \pi) \equivalent c_0  n^{-8}
e^{\alpha n}$ quand $n$ tend vers l'infini. Comme la partie imaginaire de
$\alpha$ n'est pas un multiple entier de $\pi$, il existe une suite strictement
croissante $\varphi(n)$ d'entiers tels que l'argument de 
$c_0 e^{\alpha \varphi(n)}$, vu modulo $2 \pi$, ait une limite autre que $\pm \pi /
2$. On a alors~:
$$\lim_{n \to \infty} \vert \Svetunphi(1) \vert ^{1/\varphi(n)} = 
e^{\Reelle(\alpha)}.$$
Le choix $a=20$ donne $\Reelle(\alpha) = -22,02\ldots$ d'o\`u 
$\Reelle(\alpha) + a + 2 < 0$. Donc 
la forme lin\'eaire $d_{\varphi(n)} ^{22} \Svetunphi(1)$ en 
1, $\zeta(5)$, $\zeta(7)$, \ldots, $\zeta(21)$, \`a coefficients entiers, tend
vers 0 quand $n$ tend vers l'infini et est non nulle pour $n$ assez grand. Cela
montre que l'un au moins parmi $\zeta(5)$, $\zeta(7)$, \ldots, $\zeta(21)$
est irrationnel.

\begin{rema} Si on savait d\'emontrer la conjecture mentionn\'ee 
ci-dessus (i.e.
que $2 d_n ^{a+1} \Pvetun_j(1)$ est un entier pour tout $j$), on pourrait 
(\cite{theseTanguy}, \S 5.1) appliquer
la m\^eme m\'ethode avec $a=18$, et d\'emontrer ainsi que 
l'un au moins des nombres $\zeta(5)$, $\zeta(7)$, \ldots, $\zeta(19)$,
est irrationnel.
\end{rema}


\mathversion{bold}
\subsection{Irrationalit\'e d'un nombre parmi $\zeta(5)$, $\zeta(7)$, 
$\zeta(9)$ et $\zeta(11)$}
\mathversion{normal}

La structure de la preuve est la m\^eme que dans le paragraphe pr\'ec\'edent.
La diff\'erence principale vient de d\'enominateurs nettement plus petits,
gr\^ace \`a une \'etude fine de leurs valuations $p$-adiques et \`a
l'utilisation d'une fraction rationnelle modifi\'ee~:
$$\Rzu(k) = \frac{\prod_{u=1} ^{10} ((13+2u)n)!}{(27n)!^6}(37n+2k)
\frac{(k-27n)_{27n} ^3 (k+37n+1)_{27n} ^3 }{\prod_{u=1} ^{10}
(k+(12-u)n)_{(13+2u)n+1}}.$$
Pour $\vert z \vert \geq 1$ on pose
$\Szu(z) = \frac12 \sum_{k=1} ^{\infty} \Rzu  '' (k) z^{-k}$. 
La d\'ecomposition  en \'el\'ements simples
$\Rzu(k) = \sum_{j=1} ^{10} \sum_{i=(j+1)n} ^{(36-j)n}
\frac{\cijzu}{(k+i)^{j}}$ d\'efinit les $\cijzu$ \`a partir desquels on construit
les polyn\^omes 
$\Pzu_j(z) =  \sum_{i = (j+1)n} ^{(36-j)n} \cijzu z^{i}$
pour $j \in \{1, 2, \ldots, 10\}$ et
$$\Pzu_0(z) = - \sum_{\ell=0} ^{35n-1} \left( \sum_{j=1} ^{10} \, \,
\sum_{i=\max((j+1)n,\ell+1)} ^{(36-j)n}
\frac{j(j+1)\cijzu}{2(i-\ell)^{j+2}} \right) z^\ell.$$ 
On a alors
$ \Szu(z) = \Pzu_0(z) + \sum_{j=1} ^{10} \frac{j(j+1)}{2} 
\Pzu_j(z) \Li_{j+2} (1/z)$.

Le probl\`eme est de majorer de
fa\c{c}on tr\`es pr\'ecise le d\'enominateur des rationnels $\cijzu$.
En suivant la m\'ethode 
utilis\'ee pour d\'emontrer le lemme~\ref{lemmeNikishin}, on obtiendrait 
$d_{33n}^{10-j} \cijzu \in \Z$ pour tous $i$ et
$j$. Une \'etude fine de la valuation $p$-adique des coefficients binomiaux
permet d'obtenir un d\'enominateur nettement plus petit~: on trouve un
entier $\Phi_n$ ``assez grand''
tel que $d_{33n}^{10-j} \Phi_n ^{-1} \cijzu \in \Z$. On en
d\'eduit directement que
$2d_{35n} ^3 d_{34n} d_{33n} ^8 \Phi_n^{-1} \Pzu_j(z)$ est \`a coefficients
entiers pour tout $j \in \{ 0, 1, \ldots, 10\}$.


La sym\'etrie $\Rzu(-37n-k) = - \Rzu(k)$ 
donne $z^{37n} \Pzu_j(1/z) = (-1)^{j+1} \Pzu_j(z)$, d'o\`u 
$\Pzu_j(1)=0$ pour $j = 2, 4, \ldots, 10$. 
En outre on a  $\Pzu_1(1)=0$ car $\Rzu(k) = \gdo(k^{-2})$ quand $k$ tend vers
l'infini. 
Donc $\Szu(1)$ est une forme lin\'eaire en 
1, $\zeta(5)$, $\zeta(7)$, $\zeta(9)$ et $\zeta(11)$.
Pour l'estimer, et d\'emontrer qu'elle est non nulle pour une infinit\'e de $n$,
on transforme $\Szu(1)$ en une int\'egrale complexe, \`a laquelle on
applique la m\'ethode du col (voir \cite{Zudilincentqc}, \S 2). On obtient
les comportements asymptotiques suivants quand $n$ tend vers
l'infini~:
$\limsup \vert \Szu(1) \vert ^{1/n} \leq e^{-227,58...}$, 
$\limsup \vert \Phi_n ^{-1} \vert ^{1/n} \leq e^{-176,75...}$ et 
$(d_{35n} ^3 d_{34n} d_{33n} ^8)
^{1/n} \to e^{403}$. Comme $403 < 227,58 +176,75$ on obtient la conclusion
cherch\'ee.

\begin{rema} \label{remaneuf}
Zudilin conjecture (\cite{Zudilincinqaout}, \S 9) que des
compensations ont lieu quand $z=1$, ce qui permettrait
de trouver un d\'enominateur plus petit pour les $P_j(1)$.
Peut-\^etre pourrait-on alors d\'emontrer que parmi $\zeta(5)$, $\zeta(7)$
et $\zeta(9)$ l'un au moins est irrationnel.
\end{rema}

\begin{rema} En utilisant des m\'ethodes similaires, on peut d\'emontrer 
\cite{Zudilincentqc} que pour tout $\ell \geq 1$ impair l'un au moins des
nombres $\zeta(\ell +2)$, $\zeta(\ell+4)$, \ldots, $\zeta(8 \ell -1)$, est
irrationnel.
\end{rema}

\hspace{-\parindent}St\'{e}phane Fischler 
 
\hspace{-\parindent}D\'{e}partement de Math\'{e}matiques et Applications
 
\hspace{-\parindent}\'{E}cole Normale Sup\'{e}rieure 
 
\hspace{-\parindent}45, rue d'Ulm 
 
\hspace{-\parindent}75230 Paris Cedex 05, France 
  
\hspace{-\parindent}fischler@dma.ens.fr 

\hspace{-\parindent}http://www.dma.ens.fr/$\sim$fischler/

\end{document}